\documentclass{CSML}
\pdfoutput=1

\usepackage{lastpage}
\newcommand{\dOi}{14(2:13)2018}
\lmcsheading{}{1--\pageref{LastPage}}{}{}%
{Apr.~12, 2017}{May~22, 2018}{}

\usepackage{hyperref}
\hypersetup{hidelinks}
\usepackage{array}
\usepackage{amsthm,amsmath,amssymb}
\usepackage{mathrsfs,xcolor}
\usepackage{cite}

\newtheorem*{lemma-canonisation}{Lemma~\ref{lem:canonization}}
\newtheorem*{lemma-strong-homo-small}{Lemma~\ref{lem:strong-homomorphism-from-small-substructures}}
\newtheorem*{proposition-fin-bound}{Proposition~\ref{prop:fin-bound}}
\newtheorem*{thm-our-tractability}{Theorem~\ref{thm:our-tractability}}

\newcommand{\ignore}[1]{}
\DeclareMathOperator{\Aut}{Aut}
\DeclareMathOperator{\Clo}{Clo}
\DeclareMathOperator{\Pol}{Pol}
\DeclareMathOperator{\End}{End}
\DeclareMathOperator{\Csp}{CSP}

\DeclareMathOperator{\Sym}{Sym}
\DeclareMathOperator{\HS}{\mathsf{HS}}
\DeclareMathOperator{\HSPfin}{\mathsf{HSP}_{fin}}
\DeclareMathOperator{\Ima}{Im}
\renewcommand{\Im}{\Ima} 
\DeclareMathOperator{\Comp}{Comp}

\newcommand\TB[2][B]{T_{\mathbb{#1},#2}}
\newcommand\orbitclone[2][C]{\mathscr{#1}^\typ_{#2}}

\newcommand\tuple[1]{\bar #1}
\newcommand{\algebra}[1]{{\underline #1}}
\newcommand{\algA}{\algebra A}
\newcommand{\algB}{\algebra B}

\newcommand{\algS}{\algebra S}
\newcommand{\scrC}{\mathscr C}
\newcommand{\scrD}{\mathscr D}

\newcommand{\mA}{\mathbb A}
\newcommand{\mB}{\mathbb B}
\newcommand{\mC}{\mathbb C}
\newcommand{\mD}{\mathbb D}
\newcommand{\mN}{\mathbb N}
\newcommand\Projs{\mathscr P}
\newcommand\typ{\mathrm{typ}}

\newcommand\tp{\mathrm{tp}}

\newcommand\red[1]{#1}
\definecolor{col}{HTML}{7500FF}
\newcommand\purple[1]{#1}

\newcommand\restr[2]{#1|_{#2}}
\newcommand\cl[2]{#1/\kern-4pt\sim_{#2}}

\begin{document}

\title{A Dichotomy for First-Order Reducts of Unary Structures}
\author{Manuel Bodirsky \and Antoine Mottet}
\address{Institut für Algebra, Technische Universität Dresden}
\email{\{manuel.bodirsky,antoine.mottet\}@tu-dresden.de}
\thanks{Both authors have received funding from the ERC under the European Community's Seventh Framework Programme (Grant Agreement no. 681988, CSP-Infinity), and the DFG-funded project `Homogene Strukturen, 
Bedingungserf\"ullungsprobleme, und topologische Klone' (Project number 622397). Antoine Mottet is supported by DFG Graduiertenkolleg 1763 (QuantLA)}

\maketitle

\begin{abstract}
Many natural decision problems can be formulated as constraint satisfaction problems for reducts $\mA$ of finitely bounded homogeneous structures.
This class of problems is a large generalisation of the class of CSPs over finite domains.
Our first result is a general polynomial-time reduction from such infinite-domain CSPs to finite-domain CSPs. We use this reduction
to obtain new powerful polynomial-time tractability conditions that can
be expressed in terms of the topological polymorphism clone of $\mA$. Moreover, we
study the subclass $\mathcal C$ of CSPs for structures $\mA$ that are \purple{reducts of} a structure with a unary language.
Also this class $\mathcal C$ properly extends the class of all finite-domain CSPs.
We apply our new tractability conditions to prove the general tractability conjecture of Bodirsky and Pinsker for reducts of finitely bounded homogeneous structures
for the class $\mathcal C$.
\end{abstract}

\section{Introduction}
Many computational problems in various areas of theoretical computer science can be formulated as \emph{constraint satisfaction problems (CSPs)}. Constraint satisfaction problems where the variables take values from a finite domain are reasonably well understood with respect to their computational complexity. 
The Feder-Vardi dichotomy conjecture for finite domain CSPs, which states that every finite-domain CSP is either in P or NP-complete,  {was recently solved independently by Bulatov~\cite{BulatovFVConjecture} and Zhuk~\cite{ZhukFVConjecture}.
The two solutions in fact prove} the stronger \emph{tractability conjecture}~\cite{JBK}
which provides an effective characterisation of those finite-domain CSPs that are NP-complete, and those that are in P. 
 {These breakthroughs have been obtained using concepts and results from universal algebra.}

The universal-algebraic approach can also be applied to classify the complexity of some classes of CSPs over infinite domains. 
This approach works particularly well if the constraints can be defined over a homogeneous structure with  {a} finite relational signature.
The class of CSPs that can be formulated this way is a considerable extension of the class of CSPs over finite domains, and captures many computational problems that have been studied in various research areas. For example, almost all CSPs studied in qualitative temporal and spatial reasoning belong to this class~\cite{Bodirsky-HDR}. 
The computational complexity of the CSPs where the
constraints can be defined over 
$({\mathbb Q};<)$ has been classified~\cite{tcsps-journal}.
Also constraint languages definable over the random graph~\cite{BodPin-Schaefer-both},
the random poset~\cite{KompatscherPham} or over  homogeneous tree-like structures~\cite{Phylo-Complexity}
have been classified. These results were
obtained by using a generalisation of the universal-algebraic approach from finite-domain CSPs, and structural Ramsey theory. 
However, it would be desirable to go further and to 
reduce complexity
classification tasks for CSPs over infinite domains 
to the rather more advanced classification results that
are known for finite-domain CSPs. 

In this paper, we present a first result in this direction. 
We study a  {subclass} of \emph{first-order definable structures with atoms}~\cite{locally-finite}
which recently attracted attention in automata theory
as a natural class of infinite structures with finite descriptions; see \cite{DBLP:conf/lics/BojanczykKLT13,DBLP:journals/corr/BojanczykKL14,KlinLOT14-short}, and references therein. In the
context of constraint satisfaction problems they were first studied in~\cite{locally-finite}, where the authors focussed on structures that are additionally \emph{locally finite},
a very strong restriction that is not needed for our approach.

We now state our results in more detail. 
Let $\mathbb A$ be a structure with a finite relational signature. The \emph{constraint satisfaction problem} for $\mathbb A$, denoted by $\Csp(\mathbb A)$, is the following computational problem: 

\vspace{.1cm}
\noindent{\bf Input.} A conjunction $\phi = \phi_1 \wedge \cdots \wedge \phi_m$ of atomic formulas over the signature of ${\mathbb A}$. \\
{\bf Question.} Is $\phi$ satisfiable in ${\mathbb A}$?

\vspace{.2cm}

An example of a problem that can be formulated
in this way is the three-colouring problem, 
which can be formulated as $\Csp({\mathbb A})$
where the domain $A$ of ${\mathbb A}$ has size three
and the signature contains a single binary relation  {symbol}
that denotes the inequality relation $\neq$ on $A$. 
In order to answer whether a given finite graph $G$ is 3-colourable, each edge $\{u,v\}$ of $G$
will be represented by  {an atomic formula}  of the form $u \neq v$. 
The graph is 3-colourable if and only if the corresponding
 {conjunction of atomic formulas is satisfiable}. 

Another example is the problem to decide whether a given set of constraints of the form $x=y$ or of the form $x \neq y$ has a solution
(over any domain). This can be formulated as 
$\Csp({\mathbb A})$ where
the domain of ${\mathbb A}$ 
is any countably infinite set $A$, 
and $=$ and $\neq$ denote the usual equality and inequality relations on $A$. Clearly, this problem
cannot be formulated as $\Csp({\mathbb B})$ for a 
structure ${\mathbb B}$ over a finite domain (neither can it be formulated with a locally finite structure $\mathbb B$ as considered in~\cite{locally-finite}). 

A relational structure $\mathbb A$ is called a 
\emph{first-order reduct} of a relational structure
$\mathbb B$ \purple{if} all relations of $\mathbb A$ have
a first-order definition over $\mathbb B$. 
A relational structure $\mathbb B$ is called \emph{unary} if all relations of $\mathbb B$ are unary.
We study the computational complexity of $\Csp({\mathbb A})$ for the class of all finite-signature first-order reducts $\mathbb A$ of unary structures.
Note that without loss of generality, we can focus on the case where $\mA$ has a countably infinite domain,
as any reduct of a unary structure has the same $\Csp$ as a reduct of a countably infinite structure.
In particular, our class of CSPs properly contains the class of constraint satisfaction problems over finite domains.
The structure $(A;=,\neq)$ from above is another example of a first-order reduct of a unary structure.

One of our main results is the following.

\begin{thm}\label{thm:main}
Let $A$ be a set, and let $U_1,\dots,U_n$ be subsets of $A$.
Every CSP for a first-order reduct of $(A;U_1,\dots,U_n)$ is in P or NP-complete. 
\end{thm}

In fact, we prove a stronger result and provide an algebraic characterisation of the structures in our class that have a tractable CSP. 
This characterisation coincides with the conjectured characterisation of tractable CSPs from~\cite{BPP-projective-homomorphisms} for \emph{first-order reducts of finitely bounded
homogeneous structures} (confer Conjecture~\ref{conj:inf-tract}).
In contrast to the three classification results for infinite-domain CSPs mentioned above~\cite{Phylo-Complexity,tcsps-journal,BodPin-Schaefer-both}, which involve many combinatorial case distinctions, we reduce the combinatorial work to the situation in the finite via purely conceptual arguments, building on various recent general results about topological clones~\cite{wonderland,BPP-projective-homomorphisms,BPT-decidability-of-definability}.

\vspace{.2cm}

The class of CSPs studied here 
is a subclass of 
almost any more ambitious classification project for infinite-domain CSPs.
We give some examples. 
\begin{enumerate}
\item We have already mentioned that our structures 
are \emph{structures definable with atoms} in the sense of~\cite{locally-finite}; this class of structures 
appeared under different names in various contexts (Fraenkel-Mostowski models, nominal sets, permutation models) and has not yet been studied systematically from the CSP viewpoint. 
\item An important class of infinite structures 
that have finite representations and 
where the universal-algebraic approach can be applied
is the class of 
\emph{finitely bounded homogeneous structures},
and more generally first-order reducts of finitely bounded homogeneous structures. For details we have to refer to Section~\ref{sect:csps}.
In fact, some of our results (e.g., Theorem~\ref{thm:types-tractable-implies-structure-tractable} and Theorem~\ref{thm:canonical-datalog}) hold in this general setting. 
Note that every unary structure is homogeneous, and is finitely bounded if its signature is finite.

\item Many CSPs that are of special interest in computer science and mathematics are formulated over classical structures such as the integers, the rationals, or the reals. 
Some first partial classifications are available also for such CSPs~\cite{BodMarMot,BodMarMot-temporal,JonssonThapper15}.
A natural question here is whether one can 
classify the CSP for all structures definable
over $({\mathbb Z};+,\leq)$, that is, for fragments of Presburger arithmetic, or for all structures that are definable
over $({\mathbb Q};+,\leq,1)$ or even $({\mathbb R};+,\times)$. Note that infinitely many infinite subsets are definable in these structures,
and that these subsets can be taken to be in general position (in the sense that the Boolean algebra they generate has maximal size and consists of infinite sets).
In the first structure, we can define $n$ infinite sets with infinite intersection by taking $U_i = \{z\in\mathbb Z \mid z=0\bmod p_i\}$ where $p_i$ is the $i$th prime.
In the other two structures, we can use unions of intervals.
This implies that every structure with finite unary signature is isomorphic to a first-order reduct of the structures above.

\end{enumerate}

Theorem~\ref{thm:main} is based 
on a universal-algebraic dichotomy result 
of independent interest that generalises the cyclic term theorem
of Barto and Kozik~\cite{Cyclic} and can be stated 
without reference to the CSP and to computational
complexity (Section~\ref{sect:mashups}).

\begin{thm}\label{thm:ua-result}
Let $A$ be a set, and let $\{U_1,\dots,U_n\}$ be a partition of $A$.
Let ${\scrC}$ be a closed function clone on $A$ such that the unary operations
 in ${\scrC}$ are precisely the injective functions that preserve $U_1,\dots,U_n \in A$. 
Then exactly one of the following holds:
\begin{itemize}
\item there are finitely many elements $c_1,\dots,c_k$ and a continuous clone homomorphism from $\scrC_{c_1,\dots,c_k}$ to the clone of projections on a two-element set;
\item there are unary operations $e_1,e_2\in\scrC$, an integer $m$ and an $m$-ary operation $f \in {\mathscr C}$, such that
$$\forall x_1,\dots,x_m\in A,\;e_1f(x_1,\dots,x_m) = e_2f(x_2,\dots,x_m,x_1) \, .$$
\end{itemize} 
\end{thm}

Our strategy to prove Theorem~\ref{thm:main} is then as follows: we first reduce the task to the situation 
where the polymorphism clone of the structure $\mA$ under consideration satisfies the assumptions of Theorem~\ref{thm:ua-result} (in Section~\ref{sect:lift}).
If the first case of this theorem applies, hardness of the CSP follows from general principles~\cite{wonderland}. 
If the second case applies, we associate to $\mA$ a structure over a finite domain  {which has a polynomial-time tractable CSP}, and finally
prove that $\Csp(\mA)$ reduces to this finite-domain CSP. 

Our reduction from infinite-domain CSPs to finite-domains CSPs is very general, and is another main contribution of this paper (Section~\ref{sect:csps}).
All that is needed here is that $\mA$ is definable in a finitely bounded relational structure.
For structures $\mA$ that are definable in a finitely bounded \emph{homogeneous} relational structure, this reduction yields new powerful tractability conditions,
formulated in terms of the topological polymorphism clone of $\mA$, using
known (unconditional) tractability conditions for finite-domain constraint satisfaction (in Section~\ref{sect:tractability-conditions}).

\section{Notation}
\label{sect:notation}
We denote the set $\{1,\dots,n\}$ by $[n]$.
A \emph{signature} $\tau$ is a set of function symbols and
relation symbols, 
where each symbol is associated with a natural number, called its \emph{arity}.
A \emph{$\tau$-structure} $\mA$ is a tuple
$(A; (Z^\mA)_{Z\in\tau})$ such that:
\begin{itemize}
    \item $Z^\mA\subseteq A^k$ if $Z$ is a relation symbol of arity $k$, and
    \item $Z^\mA\colon A^k\to A$ if $Z$ is a function symbol of arity $k$.
\end{itemize}
\purple{Structures} are denoted by blackboard bold letters, while their base sets are denoted by the corresponding capital roman letter. 
Let $\mA,\mB$ be $\tau$-structures with $B\subseteq A$. $\mB$ is a \emph{substructure} of $\mA$
if:
\begin{itemize}
    \item for every $R_i\in\tau$ of arity $k$, $R_i^\mB = R_i^\mA\cap B^k$, and
    \item for every function symbol $f_i\in\tau$ of arity $k$, $f_i^\mB=\restr{f_i^\mA}{B}$.
\end{itemize}

A \emph{homomorphism} between two structures $\mA,\mB$ with the same signature $\tau$ is a function $h\colon A\to B$
such that $(a_1,\dots,a_k)\in R^\mA\Rightarrow (h(a_1),\dots,h(a_k))\in R^\mB$ for every relation symbol $R\in\tau$, and such that for every function symbol $f\in\tau$,
$h(f^\mA(a_1,\dots,a_k)) = f^\mB(h(a_1),\dots,h(a_k))$.
When $h\colon A\to A$ is a homomorphism from $\mA$ to itself we say that $h$ is an \emph{endomorphism}.
An injective function $h\colon A\to B$ is an \emph{embedding} if we have $(a_1,\dots,a_k)\in R^\mA\Leftrightarrow (h(a_1),\dots,h(a_k))\in R^\mB$.
A surjective embedding is called an \emph{isomorphism}, and an \emph{automorphism} when $\mA$ equals $\mB$.

For the following definitions $\mA$ is a \emph{relational} $\tau$-structure (that is, $\tau$ only contains relation symbols).
A function $h\colon A^k\to A$  {\emph{preserves} a relation $R\subseteq A^n$}
if for all $n$-tuples $\tuple a^1,\dots,\tuple a^k$ in $R$,
we have that $h(\tuple a^1,\dots,\tuple a^k)$ is in $R$,
where $h$ is applied componentwise.
 {A function is called a \emph{polymorphism} of $\mA$ if it preserves all the relations $R^{\mA}$, for $R\in\tau$}.
We write $\Aut(\mA)$, $\End(\mA)$, $\Pol(\mA)$ for the sets of automorphisms,
endomorphisms, and polymorphisms of $\mA$. 
A relational structure $\mB$ with the same domain as $\mA$ is called a \emph{(quantifier-free) reduct} of $\mA$ if all the relations of $\mB$ are definable by (quantifier-free) first-order formulas in $\mA$.

\section{Finitely Bounded Structures}
\label{sect:csps}
A \emph{bound} of a class $\mathcal C$ of structures 
over a fixed finite relational signature $\tau$ is a finite
structure that does not embed into a structure
from $\mathcal C$, and that is minimal with this property (with respect to embeddability). 
A class of $\tau$-structures is called \emph{finitely bounded} 
if it has finitely many bounds up to isomorphism. 
    The \emph{age} of a relational structure $\mA$ is the class of all finite structures that embed into $\mA$.
Note that the age of a structure $\mA$ is finitely bounded if and only if
it has a \emph{finite universal axiomatisation}, that is,
\purple{a universal first-order sentence $\phi$
such that a finite structure $\mB$ is in the age of $\mA$ iff $\mB \models \phi$.}
A structure is called 
\emph{finitely bounded} if its age is.
The constraint satisfaction problem of a structure might be undecidable in general;
for example, $\Csp(\mathbb Z;+,\times,\{1\})$ (where $+$ and $\times$ are ternary relations
denoting the graphs of the corresponding operations)  {is the problem of satisfiability
of diophantine equations,} which is undecidable by the Matiyasevich-Robinson-Davis-Putnam theorem). 
But the CSP of a finitely bounded relational structure is in NP~\cite[Proposition 3.2.9]{Bodirsky-HDR}).

    The \emph{quantifier-free (qf-) type} of a tuple $(b_1,\dots,b_m)$, also called an \emph{$m$-type}
    in $\mB$, is 
    the set of all quantifier-free formulas
    $\phi(z_1,\dots,z_m)$ such that
    $\mB \models \phi(b_1,\dots,b_m)$. 
   If $\mB$ has a finite relational signature then 
    there are only finitely many $m$-types in $\mB$.

    Let $m$ be a positive integer.
    We define $T_{\mB,m}(\mA)$ to be the relational structure whose domain is the 
    set of $m$-types of $\mB$
    and whose relations are as follows. 
    \begin{itemize}
        \item For each symbol $R$ of $\mA$ of arity $r$, 
        let $\chi(z_1,\dots,z_r)$ 
        be a definition of $R$ in $\mB$. For
        $i \colon [r] \to [m]$ 
        we write $\langle \chi(z_{i(1)},\dots,z_{i(r)})\rangle$ 
        for the unary relation that consists of all the types 
        that contain $\chi(z_{i(1)},\dots,z_{i(r)})$,
        and add all such relations to $\TB{m}(\mA)$\footnote{\purple{In the following, we use functions to index tuples. This notation allows us
        to avoid double-subscripting and to conveniently talk about subtuples.}}.
        \item For each $r \in [m]$ and 
        $i,j \colon [r] \to [m]$, 
        define $\Comp_{i,j}$ to be the binary relation that contains
        the pairs $(p,q)$ of $m$-types such that 
        for every quantifier-free formula $\chi(z_1,\dots,z_s)$ of $\mB$
        and $t \colon [s] \to [r]$, 
        the formula $\chi(z_{it(1)},\dots,z_{it(s)})$ is in $p$ iff $\chi(z_{jt(1)},\dots,z_{jt(s)})$ is in $q$.
    \end{itemize}
    Note that if $(a_1,\dots,a_m)$ is of type $p$ and $(b_1,\dots,b_m)$ of type $q$,
	then $\Comp_{i,j}(p,q)$ holds if and only if
    $(a_{i(1)},\dots,a_{i(r)})$ and
    $(b_{j(1)},\dots,b_{j(r)})$ have the same type
    in $\mB$. Also note that if $i \colon [m] \to [m]$ 
    is the identity map, then $\Comp_{i,i}$ denotes
    the equality relation on the domain of $\TB{m}(\mA)$.
    
    The next theorem holds for arbitrary finitely bounded structures $\mB$.
 
\begin{thm}\label{thm:types-tractable-implies-structure-tractable}
Let $\mA$ be a quantifier-free 
reduct of a finitely bounded structure 
$\mB$, and suppose that $\mA$ has a finite signature. Let $m_a$ be the maximal arity of a relation in $\mA$ or $\mB$,
        and $m_b$ be the maximal size of a bound for $\mB$.
        Let $m$ be at least $\max(m_a+1,m_b,3)$.
        Then there is a polynomial-time reduction
        from $\Csp(\mA)$ to $\Csp(\TB{m}(\mA))$. 
\end{thm}

We give in the next section a sufficient condition for the existence of a polynomial-time reduction
in the other direction, from $\Csp(\TB{m}(\mA))$ to $\Csp(\mA)$.
\begin{proof}[Proof of Theorem~\ref{thm:types-tractable-implies-structure-tractable}]
		Let 
		$\Psi$ be an instance of $\Csp(\mA)$,
		and let $V=\{x_1,\dots,x_n\}$ be the variables of $\Psi$. 
		Assume without loss of generality that
		$n \geq m$. 
		We build an instance $\Phi$ of $\Csp(\TB{m}(\mA))$ as follows.
		\begin{itemize}
			\item The variable set of $\Phi$ is the set
			$I$ of  {increasing functions}\footnote{\purple{One could take $I$ to be the set of all functions $[m]\to V$ without any change to the reduction. We choose here to only take increasing functions so that the presentation of the example below is more concise.}} from $[m]$ to $V$ (where the variables are endowed with an arbitrary linear order).
            The idea of the reduction is that 
            the variable $v \in I$ 
            of $\Phi$ represents the qf-type of $(h(v(1)),\dots,h(v(m)))$
            in a satisfying assignment $h$ for $\Psi$.
\item For each conjunct $\psi$ of $\Psi$ 
			we add unary constraints to $\Phi$ as follows. The formula $\psi$ must be of the form $R(j(1),\dots,j(r))$
			where $R$ is a relation of $\mA$ and 
			$j \colon [r] \to V$. By assumption, $R$ has a qf-definition $\chi(z_1,\dots,z_r)$ over $\mB$.
			Let $v \in I$ be such that $\Im(j) \subseteq \Im(v)$.  
             Let $U$ be the relation symbol
            of $\TB{m}(\mA)$ that denotes
            the unary relation
            $\langle \chi(z_{v^{-1}j(1)},\dots,z_{v^{-1}j(r)})\rangle$. We then add $U(v)$ to $\Phi$. 
            \item Finally, for all $u,v \in I$ 
            let $k \colon [r]  \to \Im(u) \cap \Im(v)$ be a bijection.  
		We then add the constraint $\Comp_{u^{-1}k,v^{-1}k}(u,v)$.
		\end{itemize}

         {Before proving that the given reduction indeed works, we give an illustrating example.}
\begin{exa}
Let $\mA$ be $(\mN;=,\neq)$. We illustrate the reduction with the concrete instance
	$$x_1=x_2 \wedge x_2=x_3 \wedge x_{3}=x_4 \wedge x_1\neq x_4 \, .$$
	of $\Csp(\mA)$. 
	The structure $(\mN;=,\neq)$ is a reduct of 
	the homogeneous structure with domain
	$\mN$ and the empty signature, which has
	no bounds. 
    We have in this example $m=3$.

	The structure $\TB{3}(\mA)$ has a domain 
	of size five,  {where each element corresponds to a partition of $\{z_1,z_2,z_3\}$.}
     {The structure has a} unary relation $U_1$ for 
	$\langle z_2 = z_3 \rangle$,  {containing all partitions in which $z_2$ and $z_3$ belong to the same part.}
     {Similarly, the structure has a relation} $U_2$ for $\langle z_1 = z_3 \rangle$, $U_3$ for $\langle z_1 = z_2 \rangle$,
	$V_1$ for $\langle z_2 \neq z_3 \rangle$, $V_2$ for $\langle z_1 \neq z_3 \rangle$, and $V_3$ for $\langle z_1 \neq z_2 \rangle$. 
	The instance $\Phi$ of $\Csp(\TB{3}(\mA))$ 
	that our reduction creates has four variables,
	for the four order-preserving injections from
	$[3] \to \{x_1,x_2,x_3,x_4\}$ (where we order $x_1,\dots,x_4$ according to their index). Call $v_1,v_2,v_3,v_4$ these variables,
	where $\Im(v_j) = \{x_1,\dots,x_4\} \setminus \{x_j\}$. 
	We then have the following constraints in $\Phi$:
	\begin{itemize}
	\item $U_3(v_3)$ and $U_3(v_4)$ for the constraint $x_1 = x_2$ in $\Psi$;
	\item $U_1(v_4)$ and $U_3(v_1)$ for the constraint
	$x_2 = x_3$ in $\Psi$;
	\item $U_1(v_2)$ and $U_1(v_1)$ for the constraint $x_3 = x_4$ in $\Psi$; 
	\item $V_2(v_2)$ and $V_2(v_3)$ for the constraint $x_1 \neq x_4$ in $\Psi$. 
	\end{itemize}
	For the compatibility constraints we only give an example. Let $k,k' \colon [2] \to [4]$ be such that $k(1,2)=(1,3)$ and $k'(1,2)=(1,2)$. Then
	$\Comp_{k,k'}(v_4,v_2)$
	and $\Comp_{k',k'}(v_4,v_3)$ are in $\Phi$.  \qed
\end{exa}

     {We now prove that the reduction is correct.}
        Let $h\colon V\to B$  {be an assignment of the variables to the domain of $\mB$.} Let $\chi(z_1,\dots,z_r)$ be a qf-formula
        in the language of $\mB$, let $j\colon[r]\to V$,
        and let $v$ in $I$ be such that $\Im(j)\subseteq\Im(v)$.
        We first note the following property: 
        
\begin{align*}
                & \mB \models \chi(h(j(1)),\dots,h(j(r))) \\
                \text{iff} \; & (h(v(1)),\dots,h(v(m))) 
                \text{ satisfies } \chi(z_{v^{-1}j(1)},\dots,z_{v^{-1}j(r)}) \text{ in } \mB. && (\ddag)
            \end{align*}
         The property $(\ddag)$ holds since in the type of the tuple $(h(v(1)),\dots,h(v(m)))$,
        the variable $z_i$ represents the element $h(v(i))$, and therefore $z_{v^{-1}j(i)}$ represents $h(j(i))$. 

		 {($\Psi$ satisfiable implies $\Phi$ satisfiable.)} Suppose that $h\colon V\to B$ 
		satisfies 
		$\Psi$ in $\mA$. To show that 
		$\Phi$ is satisfiable in $\TB{m}(\mA)$  
		define $g \colon I \to \TB{m}(\mA)$ 
		by setting $g(v)$ to be the type of $(h(v(1)),\dots,h(v(m)))$ in $\mB$, for every $v \in I$. 
		To see that all the constraints of $\Phi$
		 are satisfied
		by $g$, let 
		$U(v)$ be a constraint in $\Phi$
		that has been introduced
		for a conjunct 
		of the form $R(j(1),\dots,j(r))$ in $\Psi$, where $j\colon[r]\to V$. 
		Let $\chi(z_1,\dots,z_r)$ be a qf-formula that defines $R$ in $\mB$. Then
		\begin{align*}
		& \quad \mA \models R(h(j(1)),\dots,h(j(r))) \\
		\Rightarrow & \quad \mB \models \chi(h(j(1)),\dots,h(j(r))) \\
		\Rightarrow & \quad \chi(z_{v^{-1}j(1)},\dots,z_{v^{-1}j(m)}) \in g(v) && \text{(because of $(\ddag)$)} \\
		\Rightarrow & \quad 
		\TB{m}(\mA) \models U(g(v)). 
		\end{align*}
		
		Next, consider a 
		constraint of the form
		$\Comp_{u^{-1}k,v^{-1}k}(u,v)$ 
		in $\Phi$, and let $r := |\!\Im(k)|$.  
		Let $\chi(z_1,\dots,z_s)$ be a qf-formula in the language of $\mB$
        and let $t\colon[s]\to[r]$.
        Suppose that $\chi(z_{u^{-1}kt(1)},\dots,z_{u^{-1}kt(s)})$ is in $g(u)$.
        From $(\ddag)$ we obtain that $\mB \models \chi(h(kt(1)),\dots,h(kt(s)))$.
        Again by $(\ddag)$ we get that $\chi(z_{v^{-1}kt(1)}, \dots, z_{v^{-1}kt(s)})$
        is in $g(v)$.
		Hence,
        \begin{align*}
          \TB{m}(\mA) \models \Comp_{u^{-1}k,v^{-1}k}(g(u),g(v)).
        \end{align*}

		 {($\Phi$ satisfiable implies $\Psi$ satisfiable.)} Conversely, suppose that $\Phi$ is satisfiable 
		in $\TB{m}(\mA)$. That is, there exists
		a map $h$ from $I$ to the  
		$m$-types in $\mB$ that satisfies all conjuncts of $\Phi$.
        We show how to obtain an assignment  
        $\{x_1,\dots,x_n\} \to A$ that satisfies $\Psi$ in $\mA$.
        Define an equivalence relation $\sim$ on $V$ as follows.
        Let $x,y\in V$.
        Let $u\in I$ be such that there are $p,q\in[m]$ such that $u(p)=x$ and $u(q)=y$.
        We define $x\sim y$ if, and only if, $h(u)$ contains the formula $z_{p} = z_{q}$.
        Note that the choice of $u$ is not important: if $u',p',q'$ are such that $u'(p')=x$ and $u'(q')=y$,
        the intersection of $\Im(u)$ and $\Im(u')$ contains $\{x,y\}$.
        Let $k\colon[r]\to \Im(u)\cap\Im(u')$ be a bijection.
        By construction, the constraint $\Comp_{u^{-1}k,u'^{-1}k}(u,u')$ is satisfied by $h$,
        which by definition of the relation means that $h(u)$ contains $z_{p}=z_{q}$ iff
        $h(u')$ contains $z_{p'}=z_{q'}$.
        
        We prove that $\sim$ is an equivalence relation. Reflexivity and
        symmetry  {are clear from the definition}.
        Assume that $x \sim y$
        and $y\sim z$.
         {Let $w\in I,p,q,r$ be such that
        $w(p)=x,w(q)=y$, and $w(r)=z$,
        which is possible since $m\geq 3$.
        Since $x\sim y$, the previous paragraph implies that $h(w)$ contains
        the formula $z_p=z_q$. Similarly, since $y\sim z$,
        the formula $z_q=z_r$ is in $h(w)$.
        Since $h(w)$ is a type, transitivity of equality implies that $z_p=z_r$
        is in $h(w)$, so that $x\sim z$.}

        Define a structure $\mC$ on $V/\kern-4pt\sim$ as follows.
        For every $k$-ary relation symbol $R$ of $\mB$ and $k$ elements $[y_1],\dots,[y_k]$ of $V/\kern-4pt\sim$,
        let $w\in I,p_1,\dots,p_k\in [m]$ be such that $w(p_i)=y_i$ (such a $w$ exists since $m\geq k$).
        Add the tuple $([y_1],\dots,[y_k])$ to $R^\mC$ if and only if $h(w)$ contains the formula $R(z_{p_1},\dots,z_{p_k})$.
        As in the paragraph above, this definition does not depend on the choice of the representatives $y_1,\dots,y_k$ or  {on the choice} of $w$.
        Proving that the definition does not depend on $w$ is straightforward.
        Suppose now that $y_1\sim y_1'$, and let $w\in I$ be such that $(w(p_1),\dots,w(p_k))=(y_1,\dots,y_k)$
        and such that $h(w)$ contains $R(z_{p_1},\dots,z_{p_k})$.
        Let $w'\in I$ be such that $(w'(q),w'(p'_1),\dots,w'(p'_k))=(y'_1,y_1,y_2,\dots,y_k)$, which is possible since $m\geq k+1$.
        We prove that $h(w')$ contains $R(z_{q},z_{p'_2},\dots,z_{p'_k})$.
        Since $y\sim y'$, we have that $h(w')$ contains $z_q=z_{p'_1}$.
        Moreover, the images of $w'$ and $w$ intersect on $y_1,\dots,y_k$, and since $h$ satisfies the $\Comp$ constraints,
        we obtain that $h(w')$ contains $R(z_{p'_1},\dots,z_{p'_k})$. It follows that $h(w')$ contains $R(z_q,z_{p'_2},\dots,z_{p'_k})$.
        Therefore, the definition of $R$ in $\mC$ does not depend on the choice of the representative for the first entry of the tuple.
        By iterating this argument for each coordinate, we obtain that $R^\mC$ is well-defined.

        We claim that $\mC$ embeds into $\mB$. Otherwise, there would exist a bound $\mD$ of size $k\leq m$ for $\mB$
        such that $\mD$ embeds into $\mC$. Let $[y_1],\dots,[y_k]$ be the elements of the image of $\mD$ under this embedding.
        Since $k\leq m$, there exist $w\in I,p_1,\dots,p_k$ such that $(w(p_1),\dots,w(p_k))=(y_1,\dots,y_k)$.
        The quantifier-free type of $([y_1],\dots,[y_k])$ in $\mC$ is in $h(w)$, by the previous paragraph.
        It follows that if $(a_1,\dots,a_m)\in B^m$ is a tuple whose quantifier-free type is $h(w)$, there
        is an embedding of $\mD$ into the substructure of $\mB$ induced by $\{a_1,\dots,a_m\}$.
        This contradicts the fact that $\mD$ does not embed into $\mB$.
        
        Let $e$ be an embedding $\mC\hookrightarrow\mB$.
        For $x\in V$ define $f(x) := e([x])$.
        We claim that $f\colon\{x_1,\dots,x_n\}\to A$ is a valid assignment for $\Psi$.
        Let $R(j(1),\dots,j(r))$ be a constraint from $\Psi$, where $j\colon[r]\to V$.
        Let $v\in I$ be such that $\Im(j)\subseteq\Im(v)$, and such that
        the constraint $\langle \chi(z_{v^{-1}j(1)},\dots,z_{v^{-1}j(r)})\rangle(v)$ is in $\Phi$.
        Since $h$ satisfies this constraint, $h(v)$ contains $\chi(z_{v^{-1}j(1)},\dots,z_{v^{-1}j(r)})$.
        It follows that $\mC\models \chi([j(1)],\dots,[j(r)])$.
        Since $e$ embeds $\mC$ into $\mB$, we obtain $\mB\models \chi(f({j(1)}),\dots,f({j(r)}))$,
        whence $\mA\models R(f({j(1)}),\dots,f({j(r)}))$, as required.

		The given reduction can be
		performed in polynomial time:
		the number of variables in the new instance is 
		in $O(n^m)$, 
		and if $c$ is the number
		of constraints in $\Psi$, then the number of constraints in $\Phi$ is in $O(cn^m+n^{2m})$.
        Each of the new constraints can be constructed in constant time. 
	\end{proof}
	
	\vspace{-.2cm}
We mention that the reduction is in fact  a \emph{first-order reduction}
(see~\cite{AtseriasBulatovDawar} for a definition).
We also note that
Theorem~\ref{thm:types-tractable-implies-structure-tractable}
applies to all CSPs that can be
described in SNP (for SNP in connection to CSPs see, e.g.,~\cite{FederVardi}).

\section{New abstract tractability conditions}
\label{sect:tractability-conditions}
We first recall basics from universal
algebra that are needed to formulate
the algebraic facts for finite-domain constraint satisfaction that are relevant for the purposes of this paper, collected in Theorem~\ref{thm:finite} (Subsection~\ref{sect:finite}). We then briefly introduce fundamental concepts 
for infinite-domain constraint satisfaction (Subsection~\ref{sect:infinite}); these concepts will also be needed in the later sections.  
Finally, we 
state and prove our new tractability conditions (Subsection~\ref{sect:abstract}). 

\subsection{Finite-Domain CSPs and Universal Algebra} 
\label{sect:finite}
An \emph{algebra} is a structure whose signature contains only function symbols,
whose interpretations are then called the \emph{(fundamental) operations} of the algebra.
Functions that are obtained as compositions of fundamental operations are called the \emph{term operations}.
A substructure of an algebra is referred to as a \emph{subalgebra}. 
An \emph{idempotent algebra} $\algA$ is an algebra whose operations are \emph{idempotent}, i.e., they satisfy $f^\algA(x,\dots,x)=x$ for all $x\in A$.
A \emph{trivial} algebra is an algebra whose operations are projections.
If $\algA$ is an algebra and $n$ is a positive integer, $\algA^n$
is defined as the algebra on $A^n$ where for each $k$-ary function symbol $f$
in the signature of $\algA$, $f^{\algA^n}$ is the function
$(A^n)^k\to A^n$ obtained by applying $f^\algA$ on tuples componentwise.
Given an algebra $\algA$, we write $\HSPfin(\algA)$ for the class of algebras 
that contains an algebra $\algebra T$ iff
there is a positive integer $n$, a subalgebra $\algS$ of $\algA^n$,
and a surjective homomorphism $\algS\to\algebra T$.
The class $\HS(\algA)$ is defined similarly, where we only allow $n=1$.

\subsubsection{Clones} A set $\scrC$ of operations over a set $D$ is called a \emph{function clone} if for all $k\geq 1$
and all $1\leq i\leq k$ it contains the \emph{projections} $\pi_i^k\colon(x_1,\dots,x_k)\mapsto x_i$,
and if $\scrC$ is closed under composition of operations.
The smallest function clone on $\{0,1\}$ is denoted by $\Projs$.
 A typical function clone is the set $\Clo(\algA)$ of term operations
of an algebra $\algA$, and indeed for every function clone $\scrC$ there exists an algebra $\algA$ such that $\scrC = \Clo(\algA)$.
 A \emph{clone homomorphism} between two clones
$\scrC$ and $\scrD$ is a map $\xi\colon\scrC\to\scrD$ such that $\xi(\pi_i^k) = \pi_i^k$ and 
$$\xi(f\circ(g_1,\dots,g_k)) = \xi(f)\circ(\xi(g_1),\dots,\xi(g_k))$$
holds for all $f,g_1,\dots,g_k\in\scrC$.
The \emph{stabilizer} of a clone $\scrC$ by a constant $c$, written $\scrC_c$, is the subclone of $\scrC$
consisting of all the operations $f$ such that $f(c,\dots,c)=c$.

\subsubsection{The Finite-Domain  {Dichotomy}}
An operation of arity $k \geq 2$ is a \emph{weak near-unanimity}~\cite{MarotiMcKenzie}
if it satisfies the equations
$$f(y,x,\dots,x)=f(x,y,x,\dots,x)=\dots=f(x,\dots,x,y)$$
for all $x,y\in A$.
An operation $f$ of arity $k \geq 2$ is said to be \emph{cyclic} (see~\cite{Cyclic})
if it satisfies $f(\tuple x) = f(\sigma\tuple x)$ for every $\tuple x=(x_1,\dots,x_k)$,
where $\sigma$ maps $(x_1,\dots,x_k)$ to $(x_2,\dots,x_k,x_1)$.
A 6-ary operation $f\colon A^6\to A$ is \emph{Siggers} (see~\cite{Siggers})
if it satisfies $f(x,y,x,z,y,z) = f(y,x,z,x,z,y)$
for all $x,y,z \in A$. 
The following is a combination of several results, old and new. 

\begin{thm}\label{thm:finite}
For finite idempotent algebras $\algA$, \ref{itm:no-homo-finite}.-\ref{itm:siggers}. are equivalent. 
\begin{enumerate}
\item\label{itm:no-homo-finite} There is no clone homomorphism from
$\Clo(\algA)$ 
to $\Projs$;
\item\label{itm:projs-not-in-hsp} $\HSPfin(\algA)$ does not contain trivial $2$-element algebras~\cite{Bir-On-the-structure};
\item\label{itm:projs-not-in-hs} $\HS(\algA)$ does not contain trivial $2$-element algebras (Prop.\ 4.14 from~\cite{BulatovJeavons});
\item\label{itm:wnu} $\Clo(\algA)$ contains a weak near-unanimity operation~\cite{MarotiMcKenzie};
\item\label{itm:cyclic} $\Clo(\algA)$ contains a cyclic operation~\cite{Cyclic};
\item\label{itm:siggers} $\Clo(\algA)$ contains a Siggers operation~\cite{Siggers}. 
\end{enumerate}
\end{thm}
For not necessarily idempotent finite algebras, this theorem fails in general, but items \ref{itm:wnu}, \ref{itm:cyclic}, \ref{itm:siggers} are still equivalent (see~\cite{wonderland}). 
 {Bulatov~\cite{BulatovFVConjecture} and Zhuk~\cite{ZhukFVConjecture} independently proved that these conditions are enough to imply a complexity dichotomy for finite-domain CSPs.}

\begin{thm}\label{conj:finite-domain-conjecture}
	Let $\mA$ be a finite structure with a finite signature, and let $\algA$ be an algebra such that $\Clo({\algA})=\Pol(\mA)$.
	Then CSP$({\mA})$ is in P if $\algA$ satisfies item \ref{itm:wnu}, \ref{itm:cyclic}, or \ref{itm:siggers} in Theorem~\ref{thm:finite},
    and $\Csp(\mA)$ is NP-complete otherwise. 
\end{thm}

\subsection{Infinite-Domain CSPs and Topology}
The concepts introduced in Section~\ref{sect:finite}
are also relevant for infinite-domain CSPs; however,
to study potential tractability conjectures and analogs of Theorem~\ref{thm:finite} for algebras on infinite domains, we  {additionally} need some model-theoretic and topological definitions that we collected in this subsection.

\subsubsection{Homogeneity}
A structure $\mA$ is said to be \emph{homogeneous} if every isomorphism between finite substructures of $\mA$ can be extended to an automorphism of $\mA$.
Examples of homogeneous structures are $(\mN;=)$ and $({\mathbb Q}; <)$.
Homogeneous structures with finite relational signature and their reducts
are examples of \emph{$\omega$-categorical} structures: a structure $\mA$ is $\omega$-categorical if all countable models of the first-order theory of $\mA$ are isomorphic. By the theorem of Engeler, Svenonius, and Ryll-Nardzewski (see, e.g.,~\cite{Hodges}), $\omega$-categoricity 
of $\mA$ is equivalent to \emph{oligomorphicity}
of $\Aut(\mA)$: a permutation group on $A$ is oligomorphic if for every $m \in {\mathbb N}$ the componentwise action on $A^m$ has finitely many orbits.
 {It is known that every homogeneous structure in a finite relational signature has \emph{quantifier elimination}, that is,
every first-order formula is equivalent over that structure to a quantifier-free formula. In particular, the first-order reducts and quantifier-free reducts
of a homogeneous structure in a finite relational signature are the same. This allows us to use Theorem~\ref{thm:types-tractable-implies-structure-tractable}
also when $\mA$ is a first-order reduct of a finitely bounded homogeneous structure $\mB$.}

\label{sect:infinite}
\subsubsection{The Topology of Pointwise Convergence}
A function clone comes naturally equipped with a topology, namely the topology of pointwise convergence where the base set is equipped with the discrete topology. This topology
is characterised by the fact that a sequence $(f_i)_{i\in\omega}$ converges to $f$ if, and only if, for every finite subset $S$ of $D$ there exists
an $i_0\in\omega$ such that for all $i\geq i_0$, we have that $f_i$ and $f$ coincide on $S$.
Given a relational structure $\mA$, the set $\Pol(\mA)$ is a function clone
over $A$ which is topologically closed in the full clone over $A$ (the clone
consisting of all the operations of finite arity over $A$).

Let $f,g$ be operations over $D$ and let $U$ be a set of permutations of $D$.
When the topological closure of
$\{\alpha\circ f\circ(\beta_1,\dots,\beta_k) \mid \alpha,\beta_i\in U\}$ contains $g$,
we say that $f$ \emph{interpolates $g$ modulo $U$}.

\subsubsection{Model-Complete Cores}
\label{sect:mc-core}
To state an infinite-domain tractability conjecture,
we also need the following concepts. 
An $\omega$-categorical structure $\mA$ is
\emph{model-complete} if every self-embedding 
of $\mA$ preserves all first-order formulas,
and it is a \emph{core} if every endomorphism of
$\mA$ is a self-embedding. 
It is known that $\mA$ is a model-complete core if and only if $\Aut(\mA)$ is dense in $\End(\mA)$; see~\cite{Bodirsky-HDR}. Note that every $\omega$-categorical homogeneous structure is model-complete.

\begin{exa}
Let $U_1,\dots,U_n$ be subsets of a countably infinite set $A$ and let $\mA$ be a reduct of $(A;U_1,\dots,U_n)$ such that $\End(\mA)$ is the set of all injections that preserve $U_1,\dots,U_n$. 
The automorphism group of $\mA$
contains all the permutations of $A$ that preserve the sets $U_1,\dots,U_n$.
Hence, if $e$ is an injective function on $A$ fixing $U_1,\dots,U_n$ and $S$ is a finite subset of $A$, then there is an automorphism $\alpha$ of $\mA$
such that $\restr{\alpha}{S}=\restr{e}{S}$. Therefore, $e$ is
in the topological closure of $\Aut(\mA)$, and $\Aut(\mA)$ is dense
in $\End(\mA)$. So $\mA$ is a model-complete core. 
\end{exa}

Every $\omega$-categorical structure
 is homomorphically equivalent to a model-complete core, which is unique up to isomorphism, and again $\omega$-categorical~\cite{Cores-journal}.

We can now state the tractability conjecture 
for reducts of finitely bounded homogeneous structures from~\cite{BPP-projective-homomorphisms}.
\begin{conj}\label{conj:inf-tract}
Let $\mA$ be a finite-signature reduct of a finitely bounded homogeneous structure. Then $\Csp(\mA)$ is in P if the model-complete core of $\mA$ does not have an expansion $\mC$ by finitely many constants so that $\Pol(\mC)$ has a continuous clone homomorphism to $\Projs$. 
\end{conj}

If $\mA$ does \emph{not} satisfy the condition 
in this conjecture, then it is known that $\Csp(\mA)$ is NP-hard~\cite{BPP-projective-homomorphisms}.

\subsubsection{Siggers operations modulo unary operations}
An important question when generalising Theorem~\ref{thm:finite} 
to infinite domains is how to replace the last three items involving weak near-unanimity operations, cyclic operations, and Siggers operations. 
  Given two unary operations $e_1,e_2$,
we say that an operation is \emph{Siggers modulo $e_1,e_2$}
if for all $x,y,z$ in $A$, we have $e_1(f(x,y,x,z,y,z))=e_2 (f(y,x,z,x,z,y))$.
\emph{Weak near-unanimity operations modulo $e_1,\dots,e_k$} and
\emph{cyclic operations modulo $e_1,e_2$} are defined similarly. A recent breakthrough by Barto and Pinsker~\cite{BartoPinskerDichotomy} gives the following dichotomy.

\begin{thm}[Theorem 1.4 in~\cite{BartoPinskerDichotomy}]\label{thm:algebraic-dichotomy-bartopinsker}
	Let $\scrC$ be a clone whose invertible elements form an oligomorphic group that is dense in the unary part of $\scrC$. Then exactly one of the following is true:
	\begin{itemize}
		\item there exists a continuous clone homomorphism from $\scrC_{c_1,\dots,c_k}$ to $\Projs$ for some $c_1,\dots,c_k$,
		\item there exists a Siggers operation in $\scrC$ modulo unary operations of $\scrC$.
	\end{itemize}
\end{thm}
Hence, the tractability conjecture for reducts $\mA$ of finitely bounded homogeneous structure (Conjecture~\ref{conj:inf-tract}) has an equivalent formulation using Siggers polymorphisms modulo endomorphisms. 

We mention another related recent result which gives an equivalent
formulation of the first item of Theorem~\ref{thm:algebraic-dichotomy-bartopinsker}
if $\mathscr C$ is additionally the polymorphism clone of a reduct of a homogeneous structures with finite relational signature~\cite{BKOPP}. 
The clone $\mathscr C$ satisfies the first item 
if and only if $\mathscr C$ has a uniformly continuous
\emph{h1 clone homomorphism} to ${\mathscr P}$.
We do not need this fact here,
but we need in Section~\ref{sect:mashups} the notion of h1 clone homomorphisms (introduced in~\cite{wonderland}):
 a map $\phi\colon\scrC\to\Projs$ is called an \emph{h1 clone homomorphism} if 
\[\phi(f\circ(\pi^k_{i_1},\dots,\pi^k_{i_n})) = \phi(f)\circ(\pi^k_{i_1},\dots,\pi^k_{i_n})\]
holds for every $n$-ary $f \in \scrC$ and all projections $\pi^k_{i_1},\dots,\pi^k_{i_n}$.

\subsubsection{Canonical Functions}
Canonical functions have been an important tool 
for classifying the complexity of classes of infinite-domain 
constraint satisfaction problems; for a more detailed
exposition, we refer the reader to~\cite{canonical}. 
Here, we 
need canonical functions to formulate our new tractability conditions for reducts of finitely bounded
homogeneous structures. 
Let $f\colon A^k\to A$, and let $G$ be a  {permutation} group on $A$.
We say that $f$ is \emph{$G$-canonical}
   if for all $m\in\mN,\alpha_1,\dots,\alpha_k\in G$ and $m$-tuples $\tuple a_1,\dots,\tuple a_k$, there exists
$\beta\in G$ such that $\beta f(\alpha_1\tuple a_1,\dots,\alpha_k\tuple a_k) = f(\tuple a_1,\dots,\tuple a_k)$.
Equivalently, this means that $f$ induces an operation $\xi^\typ_m(f)$ on orbits of $m$-tuples under $G$,
by defining $\xi^\typ_m(f)(O_1,\dots,O_k)$ as the orbit of $f(\tuple a_1,\dots,\tuple a_k)$ 
where $\tuple a_i$ is any $m$-tuple in $O_i$.
 {It is clear from the definition that the functions in $G$ are $G$-canonical,
and so are the projections on $A$.}
We can also check that the composition of $G$-canonical functions is again $G$-canonical.

If \red{$\scrC$} is a clone \red{whose operations are all canonical} with respect to  {a group $G$} then we say that
\emph{\red{$\scrC$} is canonical with respect to  {$G$}}. 
If \red{$\scrD$} is a function clone that contains $G$, the set \red{$\scrC$} of operations of \red{$\scrD$} that are $G$-canonical
is again a function clone that contains $G$. We also call the clone \red{$\scrC$} the \emph{canonical subclone of \red{$\scrD$} with respect to  {$G$}}. 

For every $m\in\mN$ the set that consists of the operations $\xi^\typ_m(f)$ for $f\in\scrC$ is a function clone
\red{$\orbitclone{m}$} on the set of orbits of $m$-tuples. The application $\xi^\typ_m\colon\scrC\to\red{\orbitclone{m}}$
is easily seen to be a continuous clone homomorphism.
We will several times use the following result, proved in~\cite{BPP-projective-homomorphisms}. 

\begin{prop}[Corollary of the proof of Proposition 6.6 in~\cite{BPP-projective-homomorphisms}]
\label{prop:bpp}
Let $\Sigma$ be a  {finite} set of equations and $\mA$ an $\omega$-categorical structure.
Suppose that for every finite $F \subset A$ and for every equation $s\approx t$ in $\Sigma$ there are $\alpha\in\Aut(\mA)$ and $s,t\in\Pol(\mA)$ such that $\alpha s|_{F} = t|_{F}$.
Then $\Sigma$ is satisfiable in $\Pol(\mA)$ modulo $\overline{\Aut(\mA)}$.
\end{prop}

A corollary of Proposition~\ref{prop:bpp} is the following (which is the original Proposition 6.6 in~\cite{BPP-projective-homomorphisms}).
Let $\mA$ be a homogeneous structure in a finite relational language of maximal arity $m$,
let \red{$\scrC$} be a closed clone of functions that are canonical with respect to $\Aut(\mA)$
and such that $\Aut(\mA)$ is contained in \red{$\scrC$}, and let $\Sigma$ be a set of equations.
If $\Sigma$ is satisfiable in \red{$\orbitclone{m}$}, then $\Sigma$ is satisfiable in \red{$\scrC$} modulo $\overline{\Aut(\mA)}$.

\subsection{Abstract Tractability Conditions}
\label{sect:abstract}
It is known that the complexity of $\Csp(\mA)$ for $\omega$-categorical structures $\mA$ only
depends on the properties of the polymorphism clone of $\mA$.
These properties can be of different nature. \emph{Abstract} properties are properties that can be expressed
using only the composition symbol and quantification over the operations in the clone, e.g., ``there exists an operation $f$, such that $f\circ (\pi_2^2,\pi_1^2) = f$.''
\emph{Topological} properties are properties that can also refer to $\Pol(\mA)$ as a topological object, e.g.,
``there exists a \emph{continuous} clone homomorphism $\Pol(\mA)\to\Projs$.''
Finally, \emph{concrete properties} are properties that refer to certain concrete operations in the polymorphism clone.
This distinction reflects the distinction between abstract clones, topological clones, and function clones.

It was shown that for an $\omega$-categorical structure $\Pol(\mA)$, the complexity of $\Csp(\mA)$ only depends on topological properties of $\Pol(\mA)$~\cite{Topo-Birk}.
However, most of the known conditions that imply that $\Csp(\mA)$ is in P are \emph{concrete} conditions. 
One notable exception is tractability from \emph{quasi near unanimity polymorphisms},
that is, polymorphisms that satisfy the identity
\begin{align*}
    f(y,x,\dots,x) &= f(x,y,x,\dots,x) = \cdots \\&= f(x,\dots,x,y) = f(x,\dots,x).
\end{align*}
If $\mA$ has a quasi near unanimity polymorphism
then $\Csp(\mA)$ is in P~\cite{BodDalJournal}. 
This tractability condition is an abstract condition (it can be rewritten using only $f$, the projection operations $\pi_1^2,\pi_2^2$, and the composition symbol).
The tractability conditions that we are able to lift from the finite are all of the abstract type.

 We first need to connect 
 the canonical polymorphisms of $\mA$ with the polymorphism clone of the associated type structure $\TB{m}(\mA)$ from Section~\ref{sect:csps}.

\begin{lem}\label{lem:polymorphisms-type-algebra}
    Let $\mA$ be a reduct of a homogeneous relational structure $\mB$ and let $\scrC$ be the polymorphisms of $\mA$ that are canonical with respect to  {$\Aut(\mB)$}. 
     {For all $m\geq 1$, we have} $\scrC^\typ_m\subseteq \Pol(\TB{m}(\mA))$.
    \end{lem}
    \begin{proof}
    We have to show that $\xi^\typ_m(f) \in \Pol(\TB{m}(\mA))$ for every $f \in \scrC$.
    Let $k$ be the arity of $f$. 
    Let $\chi(z_1,\dots,z_r)$ be a qf-definition 
    of a relation of  {$\mA$}. 
    Let $i \colon [r] \to [m]$ 
    and let $p_1,\dots,p_k$ be types in the relation 
    $\langle\chi(z_{i(1)},\dots,z_{i(r)})\rangle$ of
    $\TB{m}(\mA)$.
    Let $\tuple a^1,\dots,\tuple a^k$ be $m$-tuples whose types are $p_1,\dots,p_k$ respectively.
    Since $\mB$ is homogeneous the orbits of $m$-tuples under $\Aut(\mB)$
    and the qf-types of $\mB$ are in one-to-one correspondence and $\xi^\typ_m(f)$ can be seen as an operation on $m$-types.
    We have that $\xi^\typ_m(f)(p_1,\dots,p_k)$ is the type of $f(\tuple a^1,\dots,\tuple a^k)$ in $\mB$.
    Since $f$ preserves the relation defined by $\chi(z_{i(1)},\dots,z_{i(r)})$,
    it follows that $f(\tuple a^1,\dots,\tuple a^k)$ satisfies $\chi(z_{i(1)},\dots,z_{i(r)})$,
    which means that $\chi(z_{i(1)},\dots,z_{i(r)})$ is contained in the type of this tuple. Therefore, $\xi^\typ_m(f)$
    preserves the relations of $\TB{m}(\mA)$ of the first kind. 
   
    We now prove that $\xi^\typ_m(f)$ preserves the relations of the second kind in $T_{\mB,m}(\mA)$.
    Indeed, let $(p_1,q_1),\dots,(p_k,q_k)$ be pairs of types in $\mathrm{Comp}_{i,j}$.
    Let $(\tuple a^1,\tuple b^1),\dots,
    (\tuple a^k,\tuple b^k)$ be pairs of $m$-tuples such that $\tp(\tuple a^l)=p_l$ and $\tp(\tuple b^l)=q_l$ for all $l \in [k]$. 
    As noted above, the definition of $\mathrm{Comp}_{i,j}$
    implies that the tuples $(a^l_{i(1)},\dots,a^l_{i(r)})$
    and $(b^l_{j(1)},\dots,b^l_{j(r)})$ have the same type in $\mB$ for all $l\in[k]$.
    Since $f$ is canonical, we have that $(f(a^1_{i(1)},\dots,a^k_{i(1)}),\dots,f(a^1_{i(r)},\dots,a^k_{i(r)}))$
    has the same type as $(f(b^1_{j(1)},\dots,b^k_{j(1)}),\dots,f(b^1_{j(r)},\dots,b^k_{j(r)}))$ in $\mB$. 
    This implies that $$\mathrm{Comp}_{i,j} \big (\xi^\typ_m(f)(p_1,\dots,p_k),\xi^\typ_m(f)(q_1,\dots,q_k) \big )$$ holds in $T_{\mB,m}(\mA)$, 
    which concludes the proof. 
    \end{proof}
    
Suppose that $\mB$ is homogeneous in a finite relational language, and that $\mA$ is a reduct of $\mB$.
Suppose moreover that every polymorphism of $\mA$ is canonical with respect to  {$\Aut(\mB)$}.
The lemma above implies that $\xi^\typ_m$ is a continuous homomorphism from $\Pol(\mA)$ to $\Pol(\TB{m}(\mA))$,
if $m$ is greater than the arity of the language of $\mB$.
This in turn implies that there is a polynomial-time reduction
from $\Csp(\TB{m}(\mA))$ to $\Csp(\mA)$~\cite{Topo-Birk}. This proves the following corollary.

\begin{cor}\label{cor:equivalence-csps}
	Let $\mB$ be a finitely bounded homogeneous structure in a finite relational language,
	and let $\mA$ be a first-order reduct of $\mB$.
	Let $m$ be defined as in Theorem~\ref{thm:types-tractable-implies-structure-tractable}.
	Suppose that all the polymorphisms of $\mA$ are canonical with respect to $ {\Aut(\mB)}$.
	Then $\Csp(\mA)$ and $\Csp(\TB{m}(\mA))$ are polynomial-time equivalent.
\end{cor}

If $\mA$ is a reduct of a finitely bounded homogeneous structure $\mB$, then the inclusion in Lemma~\ref{lem:polymorphisms-type-algebra}
becomes an equality, for $m$ large enough. This fact is only mentioned for completeness and not used later, so we only sketch the proof.

\begin{lem}\label{lem:polymorphisms-type-algebra-iff}
   Let $\mA$ be a reduct of a finitely bounded homogeneous structure $\mB$ 
    and let $\scrC$ be the 
    polymorphisms of $\mA$ that are canonical with respect to  {$\Aut(\mB)$}. 
    Let $m$ be larger than each bound of $\mB$ and
    strictly larger than the maximal arity of $\mA$ 
    and $\mB$. 
    Then 
    $\scrC^\typ_m = \Pol(\TB{m}(\mA))$.
\end{lem}
\begin{proof}[Proof sketch.]
   The inclusion $\scrC^\typ_m \subseteq  \Pol(\TB{m}(\mA))$ has been shown in Lemma~\ref{lem:polymorphisms-type-algebra}. 
    For the reverse inclusion, 
    we prove that for every 
    $g \in \Pol(\TB{m}(\mA))$ there exists
    an $f \in \scrC$ such that $\xi^\typ_m(f) = g$.
    Let $k$ be the arity of $g$. 
    We prove that for every subset $F$ of $A$
    there exists a function $h$ from $F^k \to A$ 
    such that for all 
    $\bar a^1,\dots,\bar a^k \in F^m$ whose
    types are $p_1,\dots,p_k$, respectively, 
    $h(\bar a^1,\dots,\bar a^k)$
    has type $g(p_1,\dots,p_k)$. 
    A standard compactness argument then 
    shows the existence of a function $f \colon A^k \to A$ such that for all 
    $\bar a^1,\dots,\bar a^k \in A^m$ whose
    types are $p_1,\dots,p_k$, respectively, 
    $f(\bar a^1,\dots,\bar a^k)$
    has type $g(p_1,\dots,p_k)$, and such a function 
    must satisfy  $\xi_m^\typ(f) = g$.
    
    Note that we can assume without loss of generality
    that $\mB$ has for each relation symbol $R$ also a relation symbol for the complement of $R^\mB$. This does not change $\scrC^\typ_m$ or $\TB{m}(\mA)$.
The existence of a function $h$ with the properties as stated above can then be
    expressed as an instance $\Psi$ of $\Csp(\mB)$
    where 
    the variable set is $F^k$ and where
    we impose constraints from $\mB$ 
    on $\bar a^1,\dots,\bar a^k$ to enforce 
    that in any solution $h$ to this instance
    the tuple $h(\bar a^1,\dots,\bar a^k)$
    satisfies $g(p_1,\dots,p_k)$.
    Let $\Phi$ be the instance of $\Csp(\TB{m}(\mB))$
    obtained from $\Psi$ under the reduction from 
    $\Csp(\mB)$ to $\Csp(\TB{m}(\mB))$ described
    in the proof of
    Theorem~\ref{thm:types-tractable-implies-structure-tractable}. 
    The variables of $\Phi$ 
    are the order-preserving injections 
    from $[m]$ to $F^k$. 
    For $v \colon [m] \to F^k$ and $i \leq k$, 
    let $p_i$ be 
    the type of $(v(1).i,\dots,v(m).i)$ in $\mB$. 
    Then the mapping $h$ that sends $v$ to $g(p_1,\dots,p_k)$, 
    for all variables $v$ of $\Phi$, 
    is a solution to $\Phi$: 
    \begin{itemize}
    \item the constraints of $\Phi$ of the form
    $\langle \chi(\dots) \rangle(v)$ have
    been introduced to translate constraints of $\Psi$,
    and it is easy to see that they are satisfied by the choice
    of these constraints of $\Psi$ and
    by the choice of $h$. 
    \item The other constraints of $\Phi$ are of the form $\Comp_{i,j}(u,v)$ where $u,v$
    are order-preserving injections 
    from $[m]$ to $F^k$. 
    Since $g$ is a $k$-ary polymorphism of 
    $\Pol(\TB{m}(\mB))$ and hence preserves
    the relations $\Comp_{i,j}$ 
    of $\TB{m}(\mB)$, it follows that $h$
    satisfies these constraints, too. 
    \end{itemize}
    By Theorem~\ref{thm:types-tractable-implies-structure-tractable}, the instance $\Psi$ 
    of $\Csp(\mB)$ is satisfiable, too. 
\end{proof}

Using Lemma~\ref{lem:polymorphisms-type-algebra}, one can derive new tractability conditions for reducts of finitely bounded homogeneous structures.
\begin{thm}\label{thm:canonical-datalog}
Let $\mA$ be a finite-signature reduct of
a finitely bounded homogeneous structure $\mB$. 
Suppose that $\mA$ has a four-ary polymorphism $f$
and a ternary polymorphism $g$ that are 
canonical with respect to  {$\Aut(\mB)$}, that are 
weak near-unanimity operations modulo 
$\overline{\Aut(\mB)}$, 
and such that there are 
operations $e_1,e_2$ in $\overline{\Aut(\mB)}$ with 
$e_1(f(y,x,x,x)) = e_2(g(y,x,x))$ for all $x,y$. 
Then $\Csp(\mA)$ is in P. 
\end{thm}
\begin{proof}
Let $m$ be as in the statement of Theorem~\ref{thm:types-tractable-implies-structure-tractable}.
By Lemma~\ref{lem:polymorphisms-type-algebra}, 
$f' := \xi_m^\typ(f)$ and $g' := \xi_m^\typ(g)$ are polymorphisms of $T_{\mB,m}({\mathbb A})$. Moreover,
$f'$ and $g'$ must be weak near-unanimity operations,
and they satisfy $f'(y,x,x,x) = g'(y,x,x)$. 
It follows from~\cite{Maltsev-Cond} in combination with~\cite{BoundedWidth} that 
$T_{\mB,m}({\mathbb A})$ is in P (it can be solved by a Datalog program). Theorem~\ref{thm:types-tractable-implies-structure-tractable}
then implies that $\Csp(\mA)$ is in P, too. 
\end{proof}
Note that since the reduction from $\Csp(\mA)$ to $\Csp(\TB{m}(\mA))$ presented in Section~\ref{sect:csps} is a first-order reduction,
it is computable in Datalog. In particular, the hypotheses of Theorem~\ref{thm:canonical-datalog} imply that $\Csp(\mA)$ is in Datalog.
This result generalises many tractability results from the literature, for instance 
\begin{itemize}
\item the polynomial-time tractable fragments of RCC-5~\cite{RCC5JD}; 
\item the two polynomial-time algorithms for
partially-ordered time from~\cite{BroxvallJonsson};
\item polynomial-time tractable equality constraints~\cite{ecsps}; 
\item all polynomial-time tractable equivalence CSPs~\cite{equiv-csps}.
\end{itemize}
In all cases, the respective structures $\mA$ have a polymorphism $f$ such that $\xi^\typ_2(f)$ is a semilattice operation~\cite{JeavonsClosure}.
Finite structures with a semilattice polymorphism also have weak near-unanimity polymorphisms $f'$ and $g'$
that satisfy $f'(y,x,x,x) = g'(y,x,x)$ (see~\cite{Maltsev-Cond}), 
and hence $\mA$ satisfies the assumptions of Theorem~\ref{thm:canonical-datalog}. 
Using the same idea as in 
Theorem~\ref{thm:types-tractable-implies-structure-tractable}, one obtains a series of new \emph{abstract} tractability conditions:
for every known abstract tractability condition for finite domain CSPs, we obtain an abstract
tractability condition for 
reducts of finitely bounded homogeneous structures $\mB$. To show this, we first observe that
 the functions on $A$ that are canonical with respect to  {$\Aut(\mA)$} can be characterised algebraically.

\begin{prop}\label{prop:can-alg}
Let $\mA$ be a homogeneous model-complete core with a finite relational language. 
Then $f \colon A^n\to A$ is canonical with
respect to $\mA$ if and only if 
for all 
$a_1,\dots,a_n \in \End(\mA)$ 
there exist 
$e_1,e_2 \in \End(\mA)$ 
 such that 
 $$e_1 \circ f \circ (a_1,\dots,a_n) =  {e_2} \circ f \, .$$
\end{prop}
\begin{proof}
The ``if'' direction is clear. In the other direction, the assumption that $f$ is canonical gives that for every finite subset of $A$,
the equation $f\circ(a_1,\dots,a_n)\approx f$ is satisfiable modulo $\Aut(\mA)$.
By Proposition~\ref{prop:bpp}, this means that this equation is satisfiable modulo $\overline{\Aut(\mA)}=\End(\mA)$,
that is, there exist $e_1,e_2\in\End(\mA)$ such that $e_1\circ f\circ (a_1,\dots,a_n)=e_2\circ f$.
\end{proof}

Proposition~\ref{prop:can-alg} shows that 
the following close relative to Theorem~\ref{thm:canonical-datalog} is an abstract tractability condition.

\begin{thm}\label{thm:abstract-canonical-datalog}
Let $\mA$ be %
a finitely bounded homogeneous model-complete core. 
Suppose that $\mA$ has a four-ary polymorphism $f$
and a ternary polymorphism $g$ that are 
canonical with respect to  {$\Aut(\mA)$}, that are 
weak near-unanimity operations modulo $\End(\mA)$
and such that there are 
operations $e_1,e_2 \in \End(\mA)$ 
with  $e_1(f(y,x,x,x)) = e_2(g(y,x,x))$ for all $x,y$. 
Then $\Csp(\mA)$ is in P. 
\end{thm}
In the same way 
as in Theorem~\ref{thm:abstract-canonical-datalog} 
every abstract 
tractability result for finite-domain CSPs 
can be lifted to an abstract tractability condition
for $\omega$-categorical CSPs.
Note that the polynomial-time tractable cases in the classification for Graph-SAT problems~\cite{BodPin-Schaefer-both} can also be explained with the help 
of Corollary~\ref{cor:tract} below,  {using the recent solution to the finite-domain tractability conjecture}.

\begin{cor}\label{cor:tract}
Let $\mA$ be a finite-signature reduct of a finitely bounded 
homogeneous structure $\mB$, and suppose that 
$\mA$ has a Siggers (or weak nu) polymorphism $f$ modulo operations from $\overline{\Aut(\mB)}$ 
such that $f$ is canonical with respect to  {$\Aut(\mB)$}. Then $\Csp(\mA)$ is in P. 
\end{cor}
\begin{proof}
Let $m$ be as in the statement of Theorem~\ref{thm:types-tractable-implies-structure-tractable}.
By Lemma~\ref{lem:polymorphisms-type-algebra},
$\xi^\typ_m(f)$ is a polymorphism of $T_{\mB,m}(\mA)$.
Since $\xi^\typ_m(f)$ is a Siggers operation,
Theorem~\ref{conj:finite-domain-conjecture} implies that $\Csp(T_{\mB,m}(\mA))$
is in P. Then 
Theorem~\ref{thm:types-tractable-implies-structure-tractable} implies that $\Csp(\mA)$ is in P. 
\end{proof}

Finally, we mention that the non-trivial polynomial-time tractable cases 
for reducts of $({\mathbb Q};<)$ provide
examples that cannot be lifted from finite-domain tractability results this way, since the respective languages do not have non-trivial \emph{canonical} polymorphisms.

\section{Extending Clone Homomorphisms}
\label{sect:mashups}

Theorem~\ref{thm:abstract-canonical-datalog} and Corollary~\ref{cor:tract} in the previous section are of the form ``if the canonical polymorphisms of $\mA$ satisfy
some nontrivial equations, then $\Csp(\mA)$ is in P.'' We can reformulate this as follows.
Let $\scrC\subseteq\Pol(\mA)$ be the clone of polymorphisms of $\mA$ that are canonical (with respect to the automorphism group of some homogeneous finitely bounded structure $\mB$).
Corollary~\ref{cor:tract} is then equivalent to the statement: if there is no continuous clone homomorphism $\scrC\to\Projs$, then $\Csp(\mA)$ is in P.
On the other hand, we know that the existence of a uniformly continuous h1 homomorphism $\Pol(\mA)\to\Projs$ implies that $\Csp(\mA)$ is NP-complete~\cite{wonderland}.
This naturally raises the question as to when the existence of a continuous clone homomorphism $\scrC\to\Projs$ implies the existence
of a uniformly continuous h1 homomorphism $\Pol(\mA)\to\Projs$.

We focus on the case where every operation $f\in\Pol(\mA)$ interpolates modulo $\Aut(\mB)$ an operation that is canonical with respect to $\Aut(\mB)$.
We call this the \emph{canonisation property}.
In this setting, there is a natural candidate for extending a clone homomorphism $\xi\colon\scrC\to\Projs$ to $\phi\colon\Pol(\mA)\to\Projs$. Indeed,
for every $f\in\Pol(\mA)$, define $\phi(f)$ to be $\xi(g)$, where $g$ is any canonical function in $\overline{\Aut(\mB)f\Aut(\mB)}$.
We prove that when this definition does not depend on the choice of $g$, then $\phi$ is indeed a uniformly continuous h1 homomorphism.
We also give a property of $\Pol(\mA)$ --the \emph{mashup property}-- that implies that the extension $\phi$ is well-defined.

\subsection{Mashups}
We start with the central definition of this section.

\begin{defi}
	Let $g,h\colon B^k\to B$, let $1\leq\ell\leq k$, and let $r,s\in B$. 
	An operation $\omega$ is an \emph{$\ell$-mashup of $g$ and $h$ over $\{r,s\}$} if the following equations hold:
	\begin{align*}\omega(r,\dots,r,s,r,\dots,r) = g(r,\dots,r,s,r,\dots,r),\\
	  \omega(s,\dots,s,r,s,\dots,s) = h(s,\dots,s,r,s,\dots,s),\end{align*}
    where the different entry in the arguments above is the $\ell$-th entry. In case $\ell=1$, we simply write \emph{mashup}.
\end{defi}

In the following, we encourage the reader to work with the case $k=2$ in mind.
In this case,  the portion of the Cayley table of $\omega$ containing $\{r,s\}$ looks like the tables in Figure~\ref{fig:cayley}.
\begin{figure}
\begin{minipage}{0.45\textwidth}
\[\begin{array}{|c|c|c|}
\hline
 w &r & s \\
 \hline
 r & * & g(r,s)\\
 s & h(s,r) & *\\
 \hline
\end{array}\]
\end{minipage}%
\begin{minipage}{0.45\textwidth}
\[\begin{array}{|c|c|c|}
\hline
 w &r & s \\
 \hline
 r & * & h(r,s)\\
 s & g(s,r) & *\\
 \hline
\end{array}\]
\end{minipage}
\caption{The possible Cayley tables of an $\ell$-mashup of $g$ and $h$ over $\{r,s\}$, for $\ell\in\{1,2\}$. The symbol $*$ indicates that the value can be anything.}
\label{fig:cayley}
\end{figure}
The motivation for this definition is the following. If $g$ and $h$ are assumed to be projections, and we know that for all $\ell$
there exists an $\ell$-mashup of $g$ and $h$, then $g$ and $h$ must be the same projection. We now show that this remains true
when $g,h$ are operations of an algebra $\algB$ such that $\HS(\algB)$ contains a trivial algebra.

\begin{lem}\label{lem:mashup-implies-same-image}
	Let $\algB$ be an algebra, and let $g^\algB,h^\algB$ be operations of $\algB$ of arity $k$.
    Suppose that for all $\ell\in\{1,\dots,k\}$ and all distinct elements $r,s\in B$,
    there exists an operation of $\algB$ that is an $\ell$-mashup of $g^\algB$ and $h^\algB$ over $\{r,s\}$.
    Then for every trivial algebra $\algebra T$ in $\HS(\algB)$, we have $g^{\algebra T}=h^{\algebra T}$.
\end{lem}
\begin{proof}
    Let $\algS$ be a subalgebra of $\algB$ and $f$ be a \red{surjective} homomorphism $\algS\to\algebra T$.
    Suppose that $g^{\algebra T}$ is the $\ell$-th projection.
    Let $r,s$ be two elements of $\algS$ which are mapped by $f$ to different elements of $\algebra T$,
    and let $\omega^{\algB}$ be an $\ell$-mashup of $g^\algB$ and $h^\algB$ over $\{r,s\}$.
    By the assumption that $g^{\algebra T}$ is the $\ell$-th projection, we have  \[g^{\algebra T}(f(r),\dots,f(r),f(s),f(r),\dots,f(r)) = f(s),\]
    so that
    \begin{align*}
    \omega^{\algebra T}(f(r),\dots,f(r),f(s),f(r),\dots,f(r)) &= f(\omega^\algB(r,\dots,r,s,r,\dots,r))\\
                    &= f(g^\algB(r,\dots,r,s,r,\dots,r)) \\
                    &= g^{\algebra T}(f(r),\dots,f(r),f(s),f(r),\dots,f(r))\\
                    &= f(s)
    \end{align*}
    which implies that $\omega^{\algebra T}$ is the $\ell$-th projection, by the fact that $f(s)\neq f(r)$.
    Whence,
    $\omega^{\algebra T}(f(s),\dots,f(s),f(r),f(s),\dots,f(s)) = f(r),$
    and since $\omega^{\algB}$ is a mashup of $g^{\algB}$ and $h^{\algB}$ over $\{r,s\}$, we obtain
    \begin{align*}
    h^{\algebra T}(f(s),\dots,f(s),f(r),f(s),\dots,f(s)) &= f(h^\algB(s,\dots,s,r,s,\dots,s))\\
                    &= f(\omega^\algB(s,\dots,s,r,s,\dots,s)) \\
                    &= \omega^{\algebra T}(f(s),\dots,f(s),f(r),f(s),\dots,f(s))\\
                    &= f(r)
    \end{align*}
    which implies that $h^{\algebra T}$ is the $\ell$-th projection and that $g^{\algebra T}=h^{\algebra T}$ holds.
\end{proof}

Let $G$ be a \purple{permutation group on} $A$. Let $A/G$ be the set of orbits of $A$ under $G$. In the following,
the algebra $\algB$ we consider is the algebra on $A/G$ whose operations are of the form $\xi^\typ_1(f)$,
where $f$ is a $G$-canonical function on $A$. In order to prove that the mashup of two operations $\xi^\typ_1(g),\xi^\typ_1(h)$ exists in this algebra,
we therefore need to prove that there exists a $G$-canonical function $\omega$ on $A$ that induces this mashup in $\algB$.
This motivates Definitions~\ref{def:mashup-property} and~\ref{def:canonisation-property} below.
If $F\subseteq A^{A^k}$ we write $GF$
for the set
$\{\alpha\circ f \mid \alpha\in G, f\in F\}$
and $FG$ for $\{f\circ(\beta_1,\dots,\beta_k)\mid \beta_1,\dots,\beta_k\in G, f\in F\}$.
If $F$ consists of a single operation $f$, we write $GfG$ instead of $G\{f\}G$.

\begin{defi}\label{def:mashup-property}
	Let $G$ be a group of permutations on $A$, and let $\scrD$ be a clone on $A$ containing $G$.
	We say that $(G,\scrD)$ has the \emph{mashup property} if the following condition holds:
    for all $f\in\scrD$,
	all $G$-canonical functions $g,h\in\overline{GfG}$, and all orbits $O_1,O_2\in A/G$,
	we have that \red{$\scrD$} contains a $G$-canonical function $\omega$ such that $\xi^\typ_1(\omega)$ is a $1$-mashup of $\xi^\typ_1(g),\xi^\typ_1(h)$ over $\{O_1,O_2\}$.
\end{defi}

Note that in the definition above, it is equivalent to ask for the existence of a $1$-mashup or for $\ell$-mashups for all $\ell$.

\begin{defi}\label{def:canonisation-property}
    Let $G$ be a group of permutations on $A$, and let $\scrD$ be a clone on $A$ containing $G$.
    We say that $(G,\scrD)$ has the \emph{canonisation property} if for every operation $f\in \scrD$, there exists in $\overline{GfG}$
    a $G$-canonical function.
\end{defi}

For the two properties above, we say that a group $G$ has the property if the pair $(G,\scrD)$ has the property for any clone containing $G$.
It is known that every extremely amenable group has the canonisation property (Theorem~1 in~\cite{canonical}).
There are also examples of non-extremely amenable groups having the canonisation
property (we give an example in Proposition~\ref{prop:canonization} in the next section).

\begin{thm}[Mashup theorem]\label{thm:mashup}
	Let $G$ be an oligomorphic permutation group on $A$, and let $\scrD$ be a closed function clone over $A$ containing $G$.
    Let $\scrC$ be the canonical subclone
    of $\scrD$ with respect to $G$. Suppose that \red{$(G,\scrD)$} has the canonisation property and the mashup property,
	and that $\orbitclone{1}$ is idempotent.
	If there exists a clone homomorphism from $\orbitclone{1}$ to $\Projs$,
    then there exists an h1 homomorphism $\phi$ from $\scrD$ to $\Projs$.
    Moreover, $\phi$ is constant on sets of the form $\overline{GfG}$, for \red{$f\in\scrD$}.
\end{thm}
\begin{proof}
	Let $\algB$ be an algebra such that $\Clo(\algB)=\orbitclone{1}$.
	The algebra $\algB$ is idempotent by hypothesis on $\orbitclone{1}$,
	and since $G$ is oligomorphic, $\algB$ is finite.
	It follows from Theorem~\ref{thm:finite} that there exists a subalgebra $\algS$ of $\algB$ and a homomorphism $\mu\colon\algS\to\algebra T$
	where $\algebra T$ is a trivial algebra with two elements $\mu(r),\mu(s)$.
	Let $\xi'\colon \orbitclone{1}\to\Projs$ be the clone homomorphism that maps an operation $g^\algB$ to $g^{\algebra T}$.
	\red{By Lemma~\ref{lem:mashup-implies-same-image},} this homomorphism has the property that for two operations $g^\algB,h^\algB$ in $\orbitclone{1}$ of arity $k$,
	if for all $\ell\in\{1,\dots,k\}$ there is an operation in $\orbitclone{1}$ which is an $\ell$-mashup of $g^\algB$ and $h^\algB$ over $\{r,s\}$,
	then $\xi'(g)=\xi'(h)$.
	By composition of $\xi^{\typ}_1$ with $\xi'$, we obtain a clone homomorphism 
	$\xi\colon\red{\scrC}\to\Projs$.
	Define the extension $\phi$ of $\xi$ 
	to the whole clone $\scrD$ by setting $\phi(f) := \xi(g)$, where $g$
	is any $G$-canonical function that is interpolated by $f$ modulo $G$ -- such a function exists by the canonisation property, and is in $\scrD$ since
    $\scrD$ is closed.
    We claim that $\phi$ is well-defined, and that it is a h1 homomorphism.
 
	{\bf $\phi$ is well defined:} let $g,h$ be canonical and interpolated by $f$ modulo $G$.
	By the mashup property, we obtain for each $\ell\in\{1,\dots,k\}$ an operation $\omega\in\scrD$
	which is canonical and such that $\xi^\typ_1(w)$ is an $\ell$-mashup of $\xi^\typ_1(g)$ and $\xi^\typ_1(h)$ over $\{r,s\}$.
    Since this holds for all $\ell \in \{1,\dots,k\}$,
	we have by Lemma~\ref{lem:mashup-implies-same-image} that $\xi'(\xi^{\typ}_1(g))=\xi'(\xi^{\typ}_1(h))$,
	i.e., $\xi(g)=\xi(h)$ and $\phi$ is well defined.

	{\bf $\phi$ is constant on $\overline{GfG}$, for $f\in\scrD$:} let $f'$ be in $\overline{GfG}$.
	Let $g$ be canonical and interpolated by $f'$ modulo $G$.
    Note that $g$ is also interpolated by $f$ modulo $G$, so that $\phi(f) = \xi(g) = \phi(f')$.
    It follows that $\phi$ is constant on $\overline{GfG}$.
	
	{\bf $\phi$ is an h1 homomorphism:} we need to prove that $$\phi(f\circ(\pi_{i_1}^m,\dots,\pi_{i_k}^m)) = \phi(f)\circ(\pi_{i_1}^m,\dots,\pi_{i_k}^m)$$
	for every $f\in\scrD$ of arity $k\geq 1$ and every $m\geq 1$.
	Let $g\colon A^k\to A$ be canonical and interpolated by $f$ modulo $G$.
	Then $g\circ(\pi_{i_1}^m,\dots,\pi_{i_k}^m)$
	is interpolated modulo $G$ by $f\circ(\pi_{i_1}^m,\dots,\pi_{i_k}^m)$.
	So 
	\begin{align}
	\phi(f\circ(\pi_{i_1}^m,\dots,\pi_{i_k}^m)) & = \xi(g\circ(\pi_{i_1}^m,\dots,\pi_{i_k}^m)) \label{eq:link-canonical-1} \\
	& = \xi(g)\circ(\pi_{i_1}^m,\dots,\pi_{i_k}^m) \label{eq:link-canonical-2} \\
	& = \phi(f)\circ(\pi_{i_1}^m,\dots,\pi_{i_k}^m) \label{eq:link-canonical-3}
	\end{align}
	where (\ref{eq:link-canonical-1}) and (\ref{eq:link-canonical-3})
	hold by definition of $\phi$, and (\ref{eq:link-canonical-2})
	holds since $\xi$ is a clone homomorphism.
\end{proof}

In the case that $G$ is dense in the unary part of $\scrD$, we obtain the following useful corollary:

\begin{cor}\label{cor:mashup-thm-for-cores}
	Let $G$ be an oligomorphic permutation group on $A$, and let $\scrD$ be a closed function clone over $A$ containing $G$.
    Let $\scrC$ be the canonical subclone of $\scrD$ with respect to $G$.
    Suppose that $(G,\scrD)$ has the canonisation property and the mashup property,
	and that $G$ is \emph{dense} in $\scrD^{(1)}$.
	If there exists a clone homomorphism from $\orbitclone{1}$ to $\Projs$,
    then there exist finitely many elements $c_1,\dots,c_k\in A$ and a continuous clone homomorphism $\scrD_{c_1,\dots,c_k}\to\Projs$.
\end{cor}
\begin{proof}
If $G$ is dense in $\scrD^{(1)}$, then $\orbitclone{1}$ is idempotent so we can apply the previous theorem and obtain a h1 homomorphism $\phi\colon\scrD\to\Projs$
that is constant on sets of the form $\overline{GfG}$ for $f\in\scrD$. In particular, if $f\in\scrD$ and $e\in\overline{G}$, we have $\phi(e\circ f)=\phi(f)$.
Theorem~\ref{thm:algebraic-dichotomy-bartopinsker} states that either there is a clone homomorphism as in the statement, or there is a Siggers operation $w$ in $\scrD$ modulo $\overline G$.
    In the second case, $\phi(w)$ would be a Siggers operation in $\Projs$: indeed, suppose that $e_1,e_2\in\overline G$ are such that
        \[ \forall x,y,z\in A,\; e_1w(x,y,x,z,y,z) = e_2w(y,x,z,x,z,y). \]
    This property can also be written in the language of clones as 
        \[ e_1\circ w \circ(\pi_1^3,\pi_2^3,\pi_1^3,\pi_3^3,\pi_2^3,\pi_3^3) = e_2\circ w\circ (\pi_2^3,\pi_1^3,\pi_3^3,\pi_1^3,\pi_3^3,\pi_2^3). \]
    Applying $\phi$ on both sides of the equation, we obtain
    \[ \phi(w)\circ (\pi_1^3,\pi_2^3,\pi_1^3,\pi_3^3,\pi_2^3,\pi_3^3)= \phi(w)\circ (\pi_2^3,\pi_1^3,\pi_3^3,\pi_1^3,\pi_3^3,\pi_2^3),\]
    as $\phi$ is an h1 homomorphism.
    Whence, $\phi(w)$ is a Siggers operation in $\Projs$, a contradiction. So the first case of the dichotomy must apply, i.e., there are constants $c_1,\dots,c_k$
    such that $\scrC_{c_1,\dots,c_k}\to\Projs$ continuously.
\end{proof}

\subsection{Disjoint Unions of Structures}

Let $G_1,\dots,G_n$ be transitive\footnote{A permutation group $G$ on $A$ is called \emph{transitive} if for all $a,b\in A$, there exists $g\in G$ such that $g\cdot a=b$.}
oligomorphic \purple{permutation groups on} $A_1,\dots,A_n$, respectively.
We define $G_1\times\dots\times G_n$ to be the permutation group on the disjoint union $\bigcup_{i=1}^n A_i$ that contains, for each $g_1\in G_1,\dots,g_n\in G_n$,
the function $(g_1,\dots,g_n)\colon a\in A_i\mapsto g_i(a)$.
The orbits of $\bigcup_{i=1}^n A_i$ under $G_1\times\dots\times G_n$ are precisely the sets $A_1,\dots,A_n$.
We prove that if $G_1\times\dots\times G_n$ has the canonisation property, then it has the mashup property.
For the rest of this section, ``$f$ is canonical'' means ``$f$ is canonical with respect to $\prod_{i=1}^n G_i$''.

\begin{defi}[Local mashup]
    Let $G$ be a \purple{permutation group} on $A$.
	Let $g,h,\omega\colon A^k\to A$, let $S\subseteq A$,
    and let $U,V$ be two orbits of elements of $A$ under $G$.
	We say that $\omega$ is an $S$-mashup of $g$ and $h$ over $\{U,V\}$
   	iff the following holds: there exist $\alpha,\beta\in G$ such that
	for all $x_1,\dots,x_k\in S$, we have
	$$\omega(x_1,\dots,x_k) = \begin{cases}
		\alpha g(x_1,\dots,x_k)   & \text{if } (x_1,\dots,x_k)\in U\times V^{k-1}\\
		\beta h(x_1,\dots,x_k) & \text{if } (x_1,\dots,x_k)\in V\times U^{k-1}
	\end{cases}$$
\end{defi}

\begin{prop}\label{prop:building-mashups}
    Let $G_1,\dots,G_n$ be transitive oligomorphic \purple{permutation groups} on pairwise disjoint sets $A_1,\dots,A_n$.
    Let $A=\bigcup A_i$. Assume that $G:=\prod_{i=1}^n G_i$ has the canonisation property.
    Let $f\colon A^k\to A$, and let $g$ and $h$ be canonical and in $\overline{GfG}$.
	Let $i,j\in\{1,\dots,n\}$.
    There exists a canonical function $\zeta$ in $\overline{GfG}$
    which is for every finite set $S\subset A$ an $S$-mashup of $g$ and $h$ over $\{A_i,A_j\}$.
\end{prop}
\begin{proof}
    We first prove that for every finite subset $S$ of $A$, there exists in $GfG$
    an operation $\omega_S$ which is an $S$-mashup of $g$ and $h$ over $\{A_i,A_j\}$.
    Let $S\subset A$ be finite.
    Since $g$ and $h$ are in $\overline{GfG}$, there exist operations
    $\alpha,\gamma,\beta_1,\delta_1,\dots,\beta_k,\delta_k$ in $G$ such that
    \begin{align*}
    \forall x_1,\dots,x_k\in S\; \big(g(x_1,\dots,x_n) & = \gamma f(\delta_1x_1,\dots,\delta_kx_k) \\
   \wedge  \quad h(x_1,\dots,x_n) & = \alpha f(\beta_1x_1,\dots,\beta_kx_k) \big) 
   \end{align*}
    Define $\omega_S\colon A^k\to A$ by 
    \[\omega_S(x_1,\dots,x_k) := f(\epsilon_1(x_1),\dots,\epsilon_k(x_k)),\]
    where 
    \begin{align*}
    \epsilon_{1}(x) = \begin{cases}
    \delta_1(x) & \text{ if } x\in A_i \\
    \beta_1(x) & \text{ if } x \not\in A_i
    \end{cases}
    \end{align*}
    and
    \begin{align*}
    \epsilon_{\ell}(x) = \begin{cases}
    \beta_{\ell}(x) & \text{ if } x\in A_i\\
    \delta_{\ell}(x) & \text{ if } x \not\in A_i\\
    \end{cases}
    \end{align*}
    for $\ell>1$.
    It is easy to check that $\epsilon_{\ell}$ is an element of $G$, for every $\ell\in\{1,\dots,k\}$.
    This immediately gives
    \[
    \omega_S(x_1,\dots,x_k) =
    \begin{cases}
        \gamma^{-1}(g(x_1,\dots,x_k)) & \tuple x\in A_i\times (A_j)^{k-1}\\
        \alpha^{-1}(h(x_1,\dots,x_k)) & \tuple x\in A_j\times (A_i)^{k-1}
    \end{cases}
    \]
	Thus, 
	$\omega_S \in GfG$ is an $S$-mashup of $g$ and $h$ over $\{A_i,A_j\}$.

    We now prove that there exists a single operation
    which is an $S$-mashup for all finite $S\subset A$.
    Let $0,1,\dots$ be an enumeration of $A$.
    For each positive integer $m$, consider the equivalence relation on functions $\{0,\dots,m\}^k\to A$
	defined by $r\sim_m s$ iff there exists $\alpha\in G$
	such that $r = \alpha \circ s$.
	For each $m\geq 0$, this relation has finite index because the action of
	$G$ on $A$ is oligomorphic.
	Consider the following forest $\mathcal F$.
	For each $m\geq 0$
	and each operation $\omega$ which is an
	$\{0,\dots,m\}$-mashup of $g$ and $h$ in $GfG$,
	the forest $\mathcal F$ contains the vertex 
	$\cl{(\restr{\omega}{\{0,\dots,m\}})}{m}$.
	For each $m\geq 1$, if $\cl{r}{m}$ is a vertex of $\mathcal F$,
	then there is an edge $\{\cl{s}{m-1}, \cl{r}{m}\}$
	where $s=\restr{r}{\{0,\dots,m-1\}}$.
	By the first paragraph, there are infinitely many vertices in $\mathcal F$.
	Since $\sim_m$ has finite index for all $m\geq 0$,
	the forest is finitely branching, and has finitely many roots.
	By K\"onig's lemma, there exists an infinite branch in $\mathcal F$,
	which we denote by $(\cl{\omega_m}{m})_{m\geq 0}$.
	
	We now construct a chain of functions $\zeta_m\colon\{0,\dots,m\}^k\to A$
	such that $\zeta_m\subset \zeta_{m+1}$ for all $m\geq 0$,
	and such that $\zeta_m$ is $\sim_m$-equivalent to $\omega_m$.
	For $m=0$, take $\zeta_0=\omega_0$.
	Suppose that $m>0$ and that $\zeta_{m-1}$ is defined.
	There is an edge between $\omega_{m-1}$ and $\omega_m$ by hypothesis
	and $\zeta_{m-1}\sim_{m-1}\omega_{m-1}$,
	which means that there is $\alpha$ in $G$
	such that $\alpha\restr{\omega_m}{\{0,\dots,m-1\}} = \zeta_{m-1}$.
	Define $\zeta_m$ to be $\alpha\omega_m$.
	We have $\zeta_{m-1} = \restr{\zeta_m}{\{0,\dots,m-1\}}$
	and $\zeta_m\sim_m\omega_m$, as required.
	Let now $\zeta = \bigcup_{m\geq 0} \zeta_m$.
	
	It remains to prove that $\zeta$ is an $S$-mashup of $g,h$ for every finite $S\subset A$.
	Let $S$ be such a finite set, and $m$ be such that $m\geq \max(S)$.
	Since $\cl{\omega_m}{m}$ is an element of $\mathcal F$, there exists $\omega\in GfG$ that is an $\{0,\dots,m\}$-mashup of $g$ and $h$
	and such that $\restr{\omega}{\{0,\dots,m\}}=\omega_m$.
	Let $U,V$ be orbits of $A$. Let $\alpha,\beta$ be the elements in $G$ witnessing that $\omega$ is an $\{0,\dots,m\}$-mashup of $g$ and $h$.
	Let $\gamma\in G$ be such that $\omega_m = \gamma\zeta_m$.
	Then we have for all $x_1,\dots,x_k\in\{0,\dots,m\}$:
	\begin{align*}
	\zeta(x_1,\dots,x_k) &= \zeta_m(x_1,\dots,x_k)\\
		&= \gamma\omega_m(x_1,\dots,x_k)\\
		&=  \gamma\omega(x_1,\dots,x_k)\\
		&= \begin{cases} \gamma\alpha g(x_1,\dots,x_k) & \text{if } (x_1,\dots,x_k)\in U\times V^{k-1}\\
		\gamma\beta h(x_1,\dots,x_k) & {if } (x_1,\dots,x_k)\in V\times U^{k-1} \end{cases}.
	\end{align*}
	Therefore $\zeta$ is an $S$-mashup of $g,h$, with $\gamma\circ\alpha$ and $\gamma\circ\beta$ as witnesses.
	
	    Let $\zeta'$ be canonical and in $\overline{G\zeta G}$, which exists by the canonisation property for $G$.
    It is immediate that $\zeta'$ is an $S$-mashup of $g,h$ for every finite $S$.
    Moreover, $\zeta'$ is in $\overline{G\zeta G}$ and we have
    $$\overline{G\zeta G} \subseteq \overline{G\{\omega_m\colon m\geq 0\}G}
                               \subseteq \overline{GfG},$$
    so that $\zeta'$ is in $\overline{GfG}$ as required.
\end{proof}

\begin{prop}[Building Mashups]\label{prop:M-implies-mashup}
    Let $G$ be a permutation group on $A$.
    Let $g,h$ be $G$-canonical functions of arity $k$ and let $U,V$ be two orbits of elements of $A$ under $G$.
    Suppose that $\omega$ is canonical and is an $S$-mashup of $g,h$ over $\{U,V\}$ for every finite $S\subset A$.
    Then $\xi^\typ_1(\omega)$ is a mashup of $\xi^{\typ}_1(g)$ and $\xi^\typ_1(h)$ over $\{U,V\}$.
\end{prop}
\begin{proof}
    Let $x\in U,y\in V$.
    Then by definition $\xi^\typ_1(\omega)(U,V,\dots,V)$ is the orbit of $\omega(x,y,\dots,y)$ under $G$.
    Since $\omega$ is by assumption an $\{x,y\}$-mashup of $g$ and $h$, there exists an $\alpha \in G$ such that $\omega(x,y,\dots,y) = \alpha g(x,y,\dots,y)$.
    Hence, \[\xi^\typ_1(\omega)(U,V,\dots,V) = \xi^\typ_1(g)(U,V,\dots,V).\]
    We can prove similarly that
	\[\xi^\typ_1(\omega)(V,U,\dots,U) = \xi^\typ_1(h)(V,U,\dots,U),\]
    so that $\xi^\typ_1(\omega)$ is indeed a mashup of $\xi^\typ_1(g)$ and $\xi^\typ_1(h)$ over $\{U,V\}$.
\end{proof}

\begin{cor}\label{cor:group-has-property}
Let $G_1,\dots,G_n$ be transitive oligomorphic groups.
If $\prod_{i=1}^n G_i$ has the canonisation property, then it has the mashup property.
\end{cor}

\section{Reducts of Unary Structures}

In this section we study finite-signature reducts of unary structures, i.e.,
we study structures $\mA$ for which there exist subsets $U_1,\dots,U_n$ of the domain $A$ such that
the relations of $\mA$ are first-order definable in $(A;U_1,\dots,U_n)$.
We obtain a P/NP-complete dichotomy for the CSPs of reducts of unary structures, and the border between
tractability and intractability agrees with the conjectured border of Conjecture~\ref{conj:inf-tract}.

\begin{thm}\label{thm:our-tractability}
Let $\mA$ be a finite-signature reduct of a unary structure.
Then $\Csp(\mA)$ is in P if 
    the model-complete core $\mB$ 
    of $\mA$ 
    has a Siggers polymorphism modulo endomorphisms of $\mB$, 
    and is NP-complete otherwise.
\end{thm}

Without changing the class of structures that we are studying we can assume that $\{U_1,\dots,U_n\}$ forms a partition of $A$,
and that each $U_i$ is either infinite or a singleton $\{a\}$ for some $a\in A$.
We call such a partition a \emph{stabilised partition}.
Our claim above is then that for arbitrary subsets $U_1,\dots,U_n$ of $A$, there exists a stabilised partition $V_1,\dots,V_m$ of $A$
such that the structure $(A;U_1,\dots,U_n)$ is first-order definable in $(A;V_1,\dots,V_m)$. 

\subsection{The Case of Tame Endomorphisms}

 {We start by investigating reducts of unary structures whose endomorphisms are precisely
the injective operations that preserve the sets of the partition.
The milestone of this section is Theorem~\ref{thm:ua-result-restated},
which immediately implies Theorem~\ref{thm:ua-result} from the introduction.}

\begin{thm}\label{thm:ua-result-restated}
        Let $\{U_1,\dots,U_n\}$ be a stabilised partition of $A$.
Let $\mA$ be a reduct of $(A;U_1,\dots,U_n)$ such that $\End(\mA)$ is the set of injective operations that preserve $U_1,\dots,U_n$.
Let $ {\scrC}$ be the clone of polymorphisms of $\mA$ that are canonical with respect to
 {$\Aut(A;U_1,\dots,U_n)$}. Then the following are equivalent. 
\begin{enumerate}
\item\label{itm:no-homo} there is no continuous clone homomorphism from $\scrC$ to $\Projs$;
\item\label{itm:no-h1} for every $c_1,\dots,c_k\in A$, there is no continuous clone homomorphism from $\Pol(\mA)_{c_1,\dots,c_k}$ to $\Projs$;
\item\label{itm:pseudo-cyclic} $\mA$ has a cyclic (Siggers, weak near-unanimity) polymorphism modulo endomorphisms of $\mA$;
\item\label{itm:canonical-pseudo-cyclic} $\mA$ has a cyclic (Siggers, weak near-unanimity) polymorphism $f$ modulo endomorphisms of $\mA$
and $f$ is canonical with respect to  {$\Aut(A;U_1,\dots,U_n)$}.
\end{enumerate}
\end{thm}

 {The proof of the theorem will be given at the end of this subsection.
For now, we simply remark that the implications
$(\ref{itm:no-homo})\Rightarrow(\ref{itm:canonical-pseudo-cyclic})\Rightarrow(\ref{itm:pseudo-cyclic})\Rightarrow(\ref{itm:no-h1})$
are either trivial or immediate corollaries of statements from the literature.
We prove the implication from \ref{itm:no-h1} to \ref{itm:no-homo} by contraposition.}

Let $\{U_1,\dots,U_n\}$ be a stabilised partition of $A$, and let
$\scrD$ be a closed function clone over $A$ such that $\Aut(A;U_1,\dots, U_n)$ is dense in $\scrD^{(1)}$.
Let $\scrC$ be the canonical subclone of $\scrD$ with respect to $\Aut(A;U_1,\dots,U_n)$.
When $\Aut(A;U_1,\dots, U_n)$ is dense in $\scrD^{(1)}$,
for any $i\in\{1,\dots,n\}$
the map that takes $f\in\scrD$ to $\restr{f}{U_i}$
is well-defined and is a continuous clone homomorphism:
 {the restriction of some projection $\pi^n_i$ remains the same projection,
and $f|_{U_i}\circ (g_1|_{U_i},\dots,g_k|_{U_i}) = (f\circ(g_1,\dots,g_k))|_{U_i}$ holds};
the image of this clone homomorphism
is a function clone $\scrD_{U_i}$ 
over the  {set} $U_i$.
We show  {in the next two propositions} that one of the following holds: there exists some $i\in\{1,\dots,n\}$
such that $\scrD_{U_i}\to\Projs$, or $\scrD^\typ_1\to\Projs$, or $\scrD^\typ_2$ contains a cyclic operation.

Clearly, every permutation of $U_i$ is an operation in $\scrC_{U_i}$.
Such clones have been studied in~\cite{ecsps} in the context of constraint satisfaction
problems. In particular, the authors show the following.

\begin{thm}[Consequence of Theorem~7 in~\cite{ecsps}]\label{thm:classification-equality-languages}
	Let $\scrC$ be a closed clone over a countably infinite set $A$
	containing $\mathrm{Sym}(A)$.
	Then $\scrC$ has a continuous homomorphism to $\Projs$
	if and only if there is no constant unary 
	and no injective binary operation in $\scrC$.
\end{thm}

We say that an operation $f\colon A^k\to A$ is \emph{injective in its $i$th argument} if $f(\tuple a)\neq f(\tuple b)$ for all tuples $\tuple a,\tuple b$
with $a_i\neq b_i$.

\begin{prop}\label{prop:shape-canonical-functions}
         {Let $A$ be an infinite set and let $f\colon A^k\to A$ be a function that is canonical
with respect to  {$\Sym(A)$}.
Either $f$ is a constant function,
or there is an $i\in\{1,\dots,k\}$ such that $f$ is injective in its $i$th argument.}
\end{prop}
\begin{proof}
         {For two tuples $\tuple a,\tuple b$, let $I_{\tuple a,\tuple b}=\{j\in\{1,\dots,k\} \mid a_j\neq b_j\}$.
    By canonicity of $f$, if $\tuple a,\tuple b,\tuple c,\tuple d$ are such that $I_{\tuple a,\tuple b} = I_{\tuple c,\tuple d}$,
    then $f(\tuple a)=f(\tuple b)$ if and only if $f(\tuple c)=f(\tuple d)$.
    Suppose that the second case of the statement does not apply. That is, for all $i\in\{1,\dots,k\}$,
    there are tuples $\tuple a,\tuple b$ with $f(\tuple a)=f(\tuple b)$ and $i\in I_{\tuple a,\tuple b}$.
    We prove that for all $i\in\{1,\dots,k\}$, there are tuples $\tuple c,\tuple d$ such that $f(\tuple c)=f(\tuple d)$ and $I_{\tuple c,\tuple d}=\{i\}$.
    Let $i\in\{1,\dots,k\}$ be arbitrary. 
    Pick $\tuple a,\tuple b$ such that $f(\tuple a)=f(\tuple b)$ and such that $I_{\tuple a,\tuple b}$ is minimal with the property that $i\in I_{\tuple a,\tuple b}$.}

     {Suppose for contradiction that $|I_{\tuple a,\tuple b}| > 1$.
    Let $i'\in I_{\tuple a,\tuple b}\setminus\{i\}$. Let $c_1,\dots,c_k\in A$ be elements such that:
    \begin{itemize}
            \item $c_{i'} = a_{i'}$,
            \item $c_j = a_j=b_j$ for $j\in \{1,\dots,k\}\setminus I_{\tuple a,\tuple b}$,
        \item $c_j\not\in\{a_j,b_j\}$ for all $j\in I_{\tuple a,\tuple b}\setminus\{i'\}$.
    \end{itemize}
    Note that $I_{\tuple b,\tuple c} = I_{\tuple a,\tuple b}$, that $i\in I_{\tuple a,\tuple c}$, and that $I_{\tuple a,\tuple c}\subset I_{\tuple a,\tuple b}$.
    The first equality implies that $f(\tuple b)=f(\tuple c)$, by canonicity of $f$. Therefore, $f(\tuple a)=f(\tuple c)$.
    This contradicts the minimality assumption on $I_{\tuple a,\tuple b}$, which proves the claim.}

     {We can now prove that $f$ is constant. Let $\tuple a,\tuple b$ be arbitrary $k$-tuples.
    For $i\in\{0,\dots,k\}$, define $\tuple c^i$ to be the tuple $(b_1,\dots,b_i,a_{i+1},\dots,a_k)$.
    For all $i\in\{0,\dots,k-1\}$, we have $I_{\tuple c^i,\tuple c^{i+1}} = \{i+1\}$, so that by the claim above, we have $f(\tuple c^i)=f(\tuple c^{i+1})$.
    Note that $\tuple c^0=\tuple a$, and $\tuple c^{k}=\tuple b$, so that $f(\tuple a)=f(\tuple b)$.}
\end{proof}

\begin{prop}\label{prop:canonization}
        Let $U_1,\dots,U_n$ be  {a stabilised partition} of a set $A$ and let $G$ be the automorphism group of $(A;U_1,\dots,U_n)$.
	Every $f\colon A^k\to A$ interpolates modulo $G$ an operation $g\colon A^k\to A$
	that is canonical with respect to  {$\Aut(A;U_1,\dots,U_n)$}.
\end{prop}
\begin{proof}
	Let $\prec$ be any linear order on $A$
	such that if $u\in U_i,v\in U_j$ and $i<j$, then $u\prec v$, and such that
	$\prec$ is dense and without endpoints on $U_i$ whenever $U_i$ is infinite.
	The group $\Aut(A;U_1,\dots,U_n,\prec)$ is extremely amenable (this is a corollary
    of the fact that extreme amenability is preserved under direct products, and
    that the automorphism group of a countable dense linear order is extremely amenable~\cite{PestovTheorem}).
    It follows from Theorem 1 in~\cite{canonical}
	that there exists a $g$ which satisfies the conclusion of the lemma,
	except that $g$ is canonical with respect to  {$\Aut(A;U_1,\dots,U_n,\prec)$}.
	
	We prove that $g$ is also canonical with respect to  {$\Aut(A;U_1,\dots,U_n)$}.
     {Since the structure $(A;U_1,\dots,U_n)$ is homogeneous and $\omega$-categorical, the orbits of tuples in $(A;U_1,\dots,U_n)$
    can be defined by quantifier-free formulas without disjunctions (see~\cite[Corollary 6.4.2]{Hodges}).
    Since the signature of $(A;U_1,\dots,U_n)$ (including the equality relation) is binary,
    these quantifier-free formulas can be taken to be conjunctions of binary formulas.
    This implies that
    $g$ is canonical if and only if}
	for all \emph{pairs}
	$(a_1,b_1),\dots,(a_k,b_k),(c_1,d_1),\dots,(c_k,d_k)$
	such that $(a_j,b_j)$ is in the same orbit as $(c_j,d_j)$ for all $j\in \{1,\dots, k\}$,
	we have that $(g(a_1,\dots,a_k),g(b_1,\dots,b_k))$ and $(g(c_1,\dots,c_k),g(d_1,\dots,d_k))$
	are in the same orbit under $G$.
	Note that $G$ satisfies the following property:
	
	\begin{center}
		\begin{tabular}{lp{0.85\linewidth}} $(\dagger)$ & %
        two pairs $(a,b),(c,d)$ are in the same
	orbit under $G$ \\
    iff & $a$ and $c$ are in the same orbit,
	$b$ and $d$ are in the same orbit, and $a=b$ iff $c=d$.%
        \end{tabular}
	\end{center}
	
	 {Let $\tuple a, \tuple b, \tuple c, \tuple d\in A^k$ be such that for all $i\in\{1,\dots,k\}$,
	the pairs $(a_i,b_i)$ and $(c_i,d_i)$ are in the same orbit under $G$.}
	We first prove that $g(\tuple a)$ and $g(\tuple c)$ are in the same orbit under $G$.
	Using property $(\dagger)$ we know that $a_i$ and $c_i$ are in the same
	orbit under $G$. It follows that they are in the same
    orbit under $\Aut(A;U_1,\dots,U_n,\prec)$  {(because if $a_i$ and $c_i$ belong to one of the finite sets
    of the partition, they must be equal from the assumption that $\{U_1,\dots,U_n\}$ is a stabilised partition)}.
	Since $g$ is known to be canonical with respect to  {$\Aut(A;U_1,\dots,U_n,\prec)$},
	we have that $g(\tuple a)$ and $g(\tuple c)$ are in the same
	orbit under $\Aut(A;U_1,\dots,U_n,\prec)$, and therefore
	they are in the same orbit under $G$. Similarly
	we obtain that $g(\tuple b)$ and $g(\tuple d)$ are in the same
	orbit under $G$.	
	
    Therefore, it remains to check that
    $g(\tuple a)$ is equal to $g(\tuple b)$ iff  $g(\tuple c)$ equals $g(\tuple d)$.
    Suppose that $g(\tuple a)=g(\tuple b)$.
    Let $i\in\{1,\dots,n\}$ be such that all of $g(\tuple a),g(\tuple b),g(\tuple c),g(\tuple d)$ are in $U_i$.
    If $U_i$ is finite, we have that $g(\tuple a)=g(\tuple c)$ and $g(\tuple b)=g(\tuple d)$, so the property is true.
    Assume now that $U_i$ is infinite.
    For $j\in\{1,\dots,k\}$, let $e_j\in A$ be such that either:
    \begin{itemize}
        \item $a_j\prec b_j$ and $e_i$ is taken to be in $U_i$ and larger than $c_j$ and $d_j$,
        \item $b_j\prec a_j$ and $e_i$ is taken to be in $U_i$ and smaller than $c_j$ and $d_j$, or
        \item $a_j=b_j$ and $e_j=c_j$.
    \end{itemize}
    Note that if $a_j=b_j$ then $c_j=d_j$ so $e_j=c_j=d_j$.
    By definition, $(a_j,b_j), (c_j,e_j),$ and $(d_j,e_j)$ all are in the same orbit under
    $\Aut(A;U_1,\dots,U_n,\prec)$ for all $j\in\{1,\dots,k\}$.
    We therefore have that $(g(\tuple a),g(\tuple b)), (g(\tuple c),g(\tuple e)),$ and $(g(\tuple d),g(\tuple e))$ are in the same orbit under $\Aut(A;U_1,\dots,U_n,\prec)$,
    whence they are in the same orbit under $\Aut(A;U_1,\dots,U_n)$.
    Thus, if $g(\tuple a)=g(\tuple b)$ then $g(\tuple c)=g(\tuple e)$ and $g(\tuple d)=g(\tuple e)$. In particular, we have $g(\tuple c)=g(\tuple d)$.
\end{proof}

\begin{prop}\label{prop:cyclic}
	Let $U_1,\dots,U_n$ be a stabilised partition of $A$.
	Let $\scrD$ be a closed clone over $A$ such that $\Aut(A;U_1,\dots,U_n)$ is dense in $\scrD^{(1)}$.
	Let $\scrC$ be the subclone of $\scrD$ consisting of the canonical functions of $\scrD$.
	Suppose that neither $\orbitclone{1}$ 
	nor any $\scrD_{U_i}$
     {has} a continuous homomorphism to $\Projs$.
	Then $\orbitclone{2}$ contains a cyclic operation.
\end{prop}
\begin{proof}
     {Let $G$ be the permutation group $\Aut(A;U_1,\dots,U_n)$.}
	Consider any algebra
	$\algB$ such that $\Clo({\algB}) = \orbitclone{1}$. Note that $\algB$ is idempotent since
	$\scrD^{(1)} = \scrC^{(1)} 
	= \overline{G}$. Since 
	$\orbitclone{1}$ does not have a 
	homomorphism to $\Projs$, 
	Theorem~\ref{thm:finite} implies that 
	 there exists an operation $c \in \scrC$ of arity $k\geq 2$
	 such that $\xi^\typ_1(c)$ is cyclic in $\orbitclone{1}$.
	For every $i\in\{1,\dots,n\}$, by assumption $\scrD_{U_i}$ does not have a homomorphism to $\Projs$
	and since $\overline G = \scrC^{(1)}$ it cannot contain a unary constant function.
	By Theorem~\ref{thm:classification-equality-languages},
	there exists a binary operation in $\scrD$ that is injective when restricted to
	$U_i$ (if $U_i$ is finite it is a singleton by assumption, so such an operation also exists in this case). 
	One sees that such a binary operation generates a $k$-ary operation whose restriction to $U_i$ is again injective.
	Finally, by Proposition~\ref{prop:canonization}, this operation interpolates modulo $G$ a canonical function $g_i\in\scrC$ of arity $k$,
	which is still injective on $U_i$.

    We prove by induction on $m$, with $1\leq m\leq n$, that there exists in $\scrC$ an operation $g$ which is injective on $\bigcup_{i=1}^m (U_i)^k$,
    the case $m=1$ being dealt with by the paragraph above.
    So assume that the operation $g'$ is in $\scrC$ and is injective on $\bigcup_{i=1}^{m-1} (U_i)^k$.
    Define a new operation $g$ by 
    \[ g(x_1,\dots,x_k) := g_m(g'(\tuple x),g'(\sigma\tuple x),\dots,g'(\sigma^{k-1}\tuple x)),\]
    where $\sigma$ is the permutation $(x_1,\dots,x_k)\mapsto (x_2,\dots,x_k,x_1)$  {and $g_m$ is the $k$-ary canonical function whose existence
    is asserted in the previous paragraph}.
    Since $G$ is dense in $\scrD^{(1)}$, it is clear that if $\tuple x\in (U_i)^k$ and $\tuple y\in (U_j)^k$ for $i\neq j$, then $g(\tuple x)\neq g(\tuple y)$.
    If $\tuple x,\tuple y\in (U_i)^k$ are two different tuples with $i\leq m-1$, we have for all $j$ that $g'(\sigma^j\tuple x)\neq g'(\sigma^j\tuple y)$.
    Since $g_m$ is canonical and there is no constant operation in $\scrD$, this operation is injective in one of its arguments,
    by Proposition~\ref{prop:shape-canonical-functions}. 
    It follows that $g(\tuple x)\neq g(\tuple y)$.
    If $\tuple x,\tuple y\in (U_m)^k$, a similar argument works: since $g'$ is canonical and non-constant,
    it is injective in at least one of its arguments by Proposition~\ref{prop:shape-canonical-functions}.
    Whence, for at least one $j\in\{0,\dots,k-1\}$ we have that $g'(\sigma^j\tuple x)\neq g'(\sigma^j\tuple y)$,
    and by injectivity of $g_m$ on $(U_m)^k$, we obtain $g(\tuple x)\neq g(\tuple y)$.
    It follows that $g$ is canonical and injective on $\bigcup_{i=1}^{m} (U_i)^k$ as required.

	Define now $c' \in \scrC$ by
	\[c'(\tuple x) := g(c(\tuple x),c(\sigma\tuple x),\dots,c(\sigma^{k-1}\tuple x)),\]
     {where $g$ is the operation built in the previous paragraph.}
	We claim that $\xi^\typ_2(c')\in\orbitclone{2}$ is cyclic.
	It is trivial to check that $\xi^\typ_1(c')$ is cyclic in $\orbitclone{1}$.
	We now show that for all $k$-tuples $\tuple a,\tuple b$
	such that $c'(\tuple a)=c'(\tuple b)$ we have
	$c'(\sigma\tuple a) = c'(\sigma\tuple b)$.
	Suppose that $\tuple a$ and $\tuple b$ are given and
	map to the same point under $c'$.
        This means that
	\begin{align}
	& g(c(\tuple a),c(\sigma\tuple a),\dots,c(\sigma^{k-1}\tuple a)) \nonumber \\
	= \; & g(c(\tuple b),c(\sigma\tuple b),\dots,c(\sigma^{k-1}\tuple b))
	\label{eq:cycl}
	\end{align}
	Note that $(c(\tuple a),c(\sigma\tuple a),\dots,c(\sigma^{k-1}\tuple a))$
	and $g(c(\tuple b),c(\sigma\tuple b),\dots,c(\sigma^{k-1}\tuple b))$ are both tuples in $\bigcup_{i=1}^n (U_i)^k$.
	By injectivity of $g$ on this set, we therefore get that for all $j\in\{0,\dots,k-1\}$,
	the equality $c(\sigma^j\tuple a)=c(\sigma^j\tuple b)$ holds.
	By injecting this back into equation (\ref{eq:cycl}), we conclude that $c'(\sigma\tuple a)=c'(\sigma\tuple b)$.

	To show that $\xi^\typ_2(c')$ is cyclic, 
	let $(a_1,b_1)$, $\dots, (a_k,b_k)$ be pairs of elements of $A$.
	We have to show that 
	$(c'(\tuple a), c'(\tuple b))$ 
	and $(c'(\sigma\tuple a), c'(\sigma\tuple b))$ are in the same orbit under $G$. 
	Since $\xi^\typ_1(c')$ is cyclic,
	we already know that $c'(\tuple a)$ and $c'(\sigma\tuple a)$ are in the same orbit,
    and that $c'(\tuple b)$ and $c'(\sigma\tuple b)$ are in the same orbit.
    	Recall that $G$ satisfies the following property:
        two pairs $(a,b),(c,d)$ are in the same
	orbit under $G$ iff $a, c$ are in the same orbit,
	$b, d$ are in the same orbit and $a=b$ \text{iff} $c=d$. %
	So we only need to check that $c'(\tuple a)=c'(\tuple b)$ iff
	$c'(\sigma\tuple a) = c'(\sigma\tuple b)$.
	In the left-to-right direction, this is what we proved above.
	For the other direction, note that we can apply $k-1$ times
	the argument of the previous paragraph to obtain
	that if $c'(\sigma \tuple a) = c'(\sigma\tuple b)$, then
	$c'(\sigma^k(\tuple a)) = c'(\sigma^k(\tuple b))$,
	i.e., $c'(\tuple a) = c'(\tuple b)$.
\end{proof}

 {Recall that we want to prove the implication $(\ref{itm:no-h1})\Rightarrow(\ref{itm:no-homo})$
of Theorem~\ref{thm:ua-result-restated} by contraposition, that is, that if there is a continuous clone homomorphism from the canonical subclone $\scrC$ of $\Pol(\mA)$ to $\Projs$,
then there are constants $c_1,\dots,c_k\in A$ and a continuous clone homomorphism $\Pol(\mA)_{c_1,\dots,c_k}\to\Projs$.
The assumption implies that $\scrC^\typ_2$ does not contain a cyclic operation.
The previous proposition implies that $\Pol(\mA)_{U_i}$ has a continuous clone homomorphism to $\Projs$, or that $\scrC^\typ_1$ has a continuous homomorphism to $\Projs$.
In the first case, we immediately obtain a continuous clone homomorphism $\Pol(\mA)\to\Projs$.
In the second case, we apply Corollary~\ref{cor:mashup-thm-for-cores}. In order to do so, we need to prove that the automorphism group of a structure $(A;U_1,\dots,U_n)$
has the canonisation property and the mashup property. Note that  {the action of} $\Aut(A;U_1,\dots,U_n)$  {on $\bigcup_{i=1}^n U_i$} is isomorphic
to the intransitive action of $\prod_{i=1}^n \Sym(U_i)$ on $\bigcup_{i=1}^n U_i$: indeed, every element $\alpha\in\Aut(A;U_1,\dots,U_n)$
can be naturally seen as an $n$-tuple of permutations $(\alpha_1,\dots,\alpha_n)$ where $\alpha_i\in\Sym(U_i)$, and conversely
such a tuple of permutations can be interpreted as one single permutation on $\bigcup_{i=1}^n U_i=A$.
This is an isomorphism of abstract groups and moreover this isomorphism commutes with the actions of the two groups.
It follows from Corollary~\ref{cor:group-has-property} that if $\Aut(A;U_1,\dots,U_n)$ has the canonisation property,
then it has the mashup property.}

\begin{proof}[Proof of Theorem~\ref{thm:ua-result-restated}]
We prove the implications $\ref{itm:no-homo}\Rightarrow\ref{itm:canonical-pseudo-cyclic}\Rightarrow\ref{itm:pseudo-cyclic}\Rightarrow\ref{itm:no-h1}\Rightarrow\ref{itm:no-homo}$.
Let $\scrD$ be the clone of polymorphisms of $\mA$; so $\scrC$ is the subclone of $\scrD$ consisting of the operations that are canonical with respect to  {$\Aut(A;U_1,\dots,U_n)$}.

Suppose that \ref{itm:no-homo} holds, that is, there is no continuous clone homomorphism from $\scrC$ to $\Projs$.
It follows that there is no clone homomorphism from $\orbitclone{2}$ to $\Projs$.
By Theorem~\ref{thm:finite}, there exists a cyclic operation in $\orbitclone{2}$.
By Proposition~\ref{prop:bpp}, there exists an operation in $\scrC$ which is cyclic modulo endomorphisms of $\mA$.
This operation is a polymorphism of $\mA$ that is cyclic modulo endomorphisms and that is canonical with respect to  {$\Aut(A;U_1,\dots,U_n)$}.
This proves~\ref{itm:canonical-pseudo-cyclic}.

The implication $\ref{itm:canonical-pseudo-cyclic}\Rightarrow\ref{itm:pseudo-cyclic}$ is trivial.

Suppose now that $\scrD$ contains a polymorphism $f$ that is cyclic modulo endomorphisms, and let $c_1,\dots,c_k {\in A}$.
Proposition~5.6.9 in~\cite{Bodirsky-HDR} states that $\scrD_{c_1,\dots,c_k}$ contains an operation that is cyclic modulo unary operations.
This implies that there cannot be a clone homomorphism from $\scrD_{c_1,\dots,c_k}$
to $\Projs$, proving the implication $\ref{itm:pseudo-cyclic}\Rightarrow\ref{itm:no-h1}$.

It remains to prove that \ref{itm:no-h1} implies \ref{itm:no-homo}.
By contraposition, let us suppose that \ref{itm:no-homo} does not hold. Thus, there is a continuous clone homomorphism from $\scrC$ to $\Projs$.
By Proposition~6.7 in~\cite{BPP-projective-homomorphisms}, there exists a clone homomorphism from $\orbitclone{2}$ to $\Projs$.
By Theorem~\ref{thm:finite}, there is no cyclic operation in $\orbitclone{2}$. By Proposition~\ref{prop:cyclic}, either there
exists a continuous clone homomorphism $\scrD_{U_i}\to\Projs$ for some $i\in\{1,\dots,n\}$, or there is a clone homomorphism $\orbitclone{1}\to\Projs$.
In the first case we are done: we obtain by composing with $\scrD\to\scrD_{U_i}$ a continuous clone homomorphism $\scrD\to\Projs$, so \ref{itm:no-h1} does not hold.
 {Suppose we are in the second case. Proposition~\ref{prop:canonization} implies that $\Aut(A;U_1,\dots,U_n)$ has the canonisation property.
Note that this group is isomorphic as permutation group to the intransitive action of $\prod_{i=1}^n \Sym(U_i)$ on $\bigcup_{i=1}^n U_i$.
By Corollary~\ref{cor:group-has-property}, $\Aut(A;U_1,\dots,U_n)$ has the mashup property.}
Corollary~\ref{cor:mashup-thm-for-cores} implies that
there exist elements $c_1,\dots,c_k\in A$ and  a continuous clone homomorphism from $\scrD_{c_1,\dots,c_k}$ to $\Projs$.
This shows that \ref{itm:no-h1} does not hold in this case either, and concludes the proof of $\ref{itm:no-h1}\Rightarrow\ref{itm:no-homo}$.
\end{proof}

\subsection{The General Case}\label{sect:lift}

 {In this section we conclude the proof of the dichotomy theorem for reducts $\mA$ of unary structures $(A;U_1,\dots,U_n)$.
The previous section treated the special case where $\End(\mA)$ consists exactly of the injective operations preserving $U_1,\dots,U_n$.
In the following, we reduce the general case to this situation.}

 {The first step of the strategy for this is to show that we can assume without loss of generality that $\mA$ is a model-complete core.
Since reducts of unary structures are $\omega$-categorical, and since every $\omega$-categorical structure
has a model-complete core,
it suffices to prove that the model-complete core of a reduct of a unary structure
is again a reduct of a unary structure (Lemma~\ref{lem:mc-core} below).
The second step is to show that by adding constants in a suitable way, we obtain
a reduct of a unary structure which satisfies the hypothesis of the previous section \purple{(Proposition~\ref{prop:stabilising-cores})}.}

We can assume that the model-complete core $\mB$ of $\mA$ 
is a substructure of $\mA$.
Indeed, if $f\colon\mA\to\mB$ and $g\colon\mB\to\mA$ are homomorphisms,
then $f\circ g$ is an endomorphism of $\mB$, and is therefore an embedding
$\mB\hookrightarrow\mB$. This implies that $g$ is an embedding of $\mB$ into $\mA$.
We can then replace $\mB$ by the substructure of $\mA$ induced by $g(\mB)$.

\begin{lem}\label{lem:mc-core}
Let $\mA$ be a reduct of a unary structure, and let $\mB$ be the model-complete core of $\mA$.
Then $\mB$ is a reduct of a unary structure.
\end{lem}
\begin{proof}
Let $\mA$ be a reduct of $(A;U_1,\dots,U_n)$.
Suppose that $\mB$ is a substructure of $\mA$. 
Let $h$ be a homomorphism from $\mA$ to $\mB$.
We show that $\mB$ is a reduct of $(B; U_1\cap B,\dots, U_n\cap B)$.
To this end, we prove that every permutation of $B$ preserving the sets $U_1\cap B,\dots,U_n\cap B$ is an automorphism of $\mB$.
Let $\beta$ be such a permutation. Then $\beta$ can be extended by the identity to a permutation $\alpha$
of $A$ which preserves $U_1,\dots,U_n$, and therefore $\alpha$ is an automorphism of $\mA$.
Thus, $h\circ\beta=\restr{h\circ\alpha}{B}\colon \mB\to\mB$ is an endomorphism of $\mB$,
and so an embedding since $\mB$ is a model-complete core.
This implies that $\beta$ is an embedding, i.e., it is an automorphism of $\mB$. 
Note that $(B; U_1\cap B,\dots, U_n\cap B)$ is $\omega$-categorical.
It is a known corollary of the Ryll-Nardzewski theorem~\cite{Hodges}
that if a structure $\mC$ is $\omega$-categorical
and $\Aut(\mC)\subseteq\Aut(\mB)$ then $\mB$ is a reduct of $\mC$.
It follows that $\mB$ is a first-order 
reduct of $(B; U_1\cap B,\dots, U_1\cap B)$.
\end{proof}

It can be the case that $\End(\mA)$
contains more operations than the injections preserving $U_1,\dots,U_n$ even when $\mA$ is a reduct of $(A;U_1,\dots,U_n)$ which is a model-complete core.
An example is $(A; E,  {\neq})$ where $A=U_1\uplus U_2$ and $E = \{(x,y) \in A^2 \mid x\in U_1\Leftrightarrow y\in U_2\}$.
However, for every such reduct there are finitely many constants $c_1,\dots,c_n \in A$
such that the $(\mA,c_1,\dots,c_n)$ satisfies the condition of Theorem~\ref{thm:ua-result-restated}.

\begin{prop}\label{prop:stabilising-cores}
        Let $\mA$ be a reduct of a unary structure that is a model-complete core.
    \purple{There exist elements $c_1,\dots,c_n\in A$ and a stabilised partition $\{V_1,\dots,V_m\}$ of $A$ such that $(\mA,c_1,\dots,c_n)$ is a reduct
    of the unary structure $(A;V_1,\dots,V_m)$ and such that the endomorphisms of $(\mA,c_1,\dots,c_n)$
    are precisely the injective functions preserving $V_1,\dots,V_m$.}
\end{prop}
\begin{proof}
    \purple{Let $\{U_1,\dots,U_n\}$ be a stabilised partition of $A$ where $n$ is minimal with the property that $\mA$ is a reduct of $(A;U_1,\dots,U_n)$.}
    Up to a permutation of the blocks, we can assume that $U_1,\dots,U_r$ are the finite blocks of the partition.
    \purple{For every $i\in\{1,\dots,n\}$, let $c_i\in U_i$. We claim that
    $$\Aut(\mA,c_1,\dots,c_n)=\Aut(A;U_1,\dots,U_n,c_1,\dots,c_n).$$}
    If $r=n$, there is nothing to prove, because of the assumption that the sets $U_1,\dots,U_n$ are either singletons are infinite.
    Therefore, if $r=n$, we have $\Aut(\mA,c_1,\dots,c_n)=\Aut(A;U_1,\dots,U_n)$.

    We prove that $\Aut(\mA,c_1,\dots,c_r)$ preserves the binary relation
   	\[ E := \{(x,y) \in A^2\mid \forall i\in\{r+1,\dots,n\}, x\in U_i \Leftrightarrow y\in U_i\}. \]
    Let $\alpha$ be an automorphism of $\mA$.
    For $i,j\in\{r+1,\dots,n\}$, define $V_{ij}(\alpha)$ to be the set of elements of $U_i$ that are mapped to $U_j$ under $\alpha$.
   
   \medskip
    \underline{Claim 0:} for every $i\in\{r+1,\dots,n\}$ and every automorphism $\alpha$ of $\mA$, there exists a $j\in\{r+1,\dots,n\}$ such that $V_{ji}(\alpha)$ is infinite.
    \begin{proof}
      \pushQED{\hfill$\lozenge$}
     {Since $\alpha$ is a bijection, every element of $U_i$ has a preimage under $\alpha$. Since there are only finitely many sets in the partition,
    one of the sets $U_j$ contains infinitely many of those preimages, i.e., $V_{ji}(\alpha)$ is infinite.}
    \end{proof}

    \underline{Claim 1:} for every $i\in\{r+1,\dots,n\}$ and every automorphism $\alpha$ of $\mA$, the set $V_{ii}(\alpha)$ is either finite or $U_i$.
    \begin{proof}
    \pushQED{\hfill$\lozenge$}
    Let $\alpha$ be an automorphism of $\mA$, and suppose that $\emptyset\neq V_{ii}(\alpha)\neq U_i$.
    Since $V_{ii}(\alpha)\neq U_i$, there exists a $j\in\{r+1,\dots,n\}$ such that  {$V_{ij}(\alpha)\neq\emptyset$, which is equivalent to say that $V_{ji}(\alpha^{-1})\neq\emptyset$}.
         {Suppose for contradiction that $V_{ii}(\alpha)$ is infinite. Equivalently, $V_{ii}(\alpha^{-1})$ is infinite.}
     {We claim that for every finite subset $S$ of $A$, there exists an automorphism $\alpha'$ of $\mA$ such that $\alpha'(S)\cap U_j=\emptyset$.
    This is clear: let $\beta$ be an automorphism of $\mA$ that maps $S\cap U_i$ to $V_{ii}(\alpha^{-1})$ (which is possible since $V_{ii}(\alpha^{-1})$ is infinite)
    and one element of $S\cap U_j$ to $V_{ji}(\alpha^{-1})$.
    The automorphism $\alpha'_1:=\alpha^{-1}\circ\beta$ maps one point from $S\cap U_j$ to $U_i$, and maps all the elements of $S\cap U_i$ to $U_i$.
    Possibly, some elements in $S\cap U_k$ for $k\not\in\{i,j\}$ are mapped by $\alpha'_1$ to $U_j$.
    We repeat this procedure and obtain automorphisms $\alpha'_2,\dots,\alpha'_{m}$ with $m\leq |S|$, until $\alpha'_{m}(S)\cap U_j$ is empty.
    Using a standard compactness argument, we obtain \purple{an operation $e\in\overline{\Aut(\mA)}$} whose image does not intersect $U_j$.
    This is a contradiction to the minimality of the partition $\{U_1,\dots,U_n\}$: the structures $e(\mA)$ and $\mA$ are isomorphic,
    and the relations of $e(\mA)$ are definable in $(A\setminus U_j; U_1,\dots,U_{j-1},U_{j+1},\dots,U_n)$.}
    Therefore $V_{ii}(\alpha)$ is finite.
    \end{proof}

    \underline{Claim 2:} for every $i\in\{r+1,\dots,n\}$ and every automorphism $\alpha$ of $\mA$, the set $V_{ii}(\alpha)$ is either empty or $U_i$.
    \begin{proof}
    \pushQED{\hfill$\lozenge$}
    Suppose that for some $\alpha\in\Aut(\mA)$, the set $V_{ii}(\alpha)$ is not equal to $U_i$ and is not empty.
    We prove that for every $k\geq 1$, there exists an automorphism $\alpha_k$ of $\mA$ such that $|V_{ii}(\alpha_k)|\geq k$ and such that $\alpha_k$ does not preserve $U_i$.
    Let $k\geq 1$.
    By Claim 0, there exists a $j\in\{r+1,\dots,n\}$ such that $V_{ji}(\alpha)$ is infinite and by Claim 1, it must be the case that $j\neq i$.
    Note that $V_{ij}(\alpha^{-1})$ is infinite, and that $V_{ii}(\alpha^{-1})$ is not empty.
    Let $x_1,\dots,x_k$ be pairwise distinct elements in $V_{ij}(\alpha^{-1})$, and let $y\in V_{ii}(\alpha^{-1})$.
    Let $z$ be an element of $U_i$ such that $\alpha(z)\not\in U_i$, which exists since $V_{ii}(\alpha)\neq U_i$.
    Let $\beta$ be an automorphism of $\mA$ that maps $\alpha^{-1}(y)$ to $z$ and which leaves $\alpha^{-1}(x_1),\dots,\alpha^{-1}(x_k)$ fixed.
    Then $\alpha\circ\beta\circ\alpha^{-1}$ is an automorphism of $\mA$ such that $x_1,\dots,x_k\in V_{ii}(\alpha\circ\beta\circ\alpha^{-1})$
    and such that $(\alpha\circ\beta\circ\alpha^{-1})(y)\not\in U_i$. 
     
     For each $k\geq 1$, there exists by Claim 0 a $j\in\{r+1,\dots,n\}$ such that $V_{ji}(\alpha_k)$ is infinite.
     Since $\alpha_k$ does not preserve $U_i$ by assumption, $V_{ii}(\alpha_k)\neq U_i$.
     By Claim 1, $V_{ii}(\alpha_k)$ has to be finite, so $j$ is distinct from $i$.
     By the pigeonhole principle, there is a $j\in\{r+1,\dots,n\}$ distinct from $i$ such that $V_{ji}(\alpha_k)$ is infinite
     for infinitely many $k$.
    Therefore, using \purple{another} argument one can show that there is an endomorphism of $\mA$ in
    $\overline{\langle \Aut(A;U_1,\dots,U_n)\cup \{\alpha_k : k\geq 1\}\rangle}$
    whose image does not intersect $U_j$, which is a contradiction to the minimality of the partition $\{U_1,\dots,U_n\}$.
    Hence, $V_{ii}(\alpha)$ is either empty or $U_i$.
    \end{proof}

    \underline{Claim 3:} for every $i\in\{r+1,\dots,n\}$ and every automorphism $\alpha$ of $\mA$, there is exactly one $j\in\{r+1,\dots,n\}$ such that $V_{ij}(\alpha)$ is nonempty.
    \begin{proof}
    \pushQED{\hfill$\lozenge$}
    Suppose that $j,j'\in\{r+1,\dots,n\}$ are distinct and that $V_{ij}(\alpha)$ and $V_{ij'}(\alpha)$ are both nonempty,
    say that $\alpha(x)\in U_j$ and $\alpha(y)\in U_{j'}$.
    Since $V_{jj}(\alpha^{-1})$ is not $U_j$, it must be empty by Claim~2. Thus, there exists a $k$ distinct from $j$ such that $V_{jk}(\alpha^{-1})$ is infinite,
    which gives the existence of a $z\in U_j$ distinct from $\alpha(x)$ such that $\alpha^{-1}(z)\in U_k$.
    Let $\beta$ be an automorphism of $\mA$ that maps $x$ to $y$, and leaves $\alpha^{-1}(z)$ fixed.
    Then the map $\alpha\circ\beta\circ\alpha^{-1}$ maps $\alpha(x)\in U_j$ to $\alpha(y)\in U_{j'}$, and maps $z\in U_j$ to itself.
    Therefore, we have that $V_{jj}(\alpha\circ\beta\circ\alpha^{-1})$ is neither empty nor equal to $U_j$, a contradiction to the second claim.
    \end{proof}
    
    Therefore, the relation $E$ is preserved by $\Aut(\mA\purple{, c_1,\dots,c_r})$.
    This implies that each of $U_1,\dots,U_n$ is preserved by $\Aut(\mA,c_1,\dots,c_n)$.
    \purple{We obtain that $(\mA,c_1,\dots,c_n)$ is a reduct of $(A;U_1\setminus\{c_1\},\dots,U_n\setminus\{c_n\},\{c_1\},\dots,\{c_n\})$
    whose endomorphisms are precisely the injective functions that preserve this stabilised partition.}
\end{proof}

\begin{cor}\label{cor:Siggers-iff-no-homo} 
    Let $\mA$ be a reduct of a unary structure. Then there
    exists an expansion $\mC$ of the model-complete core of $\mA$
       by finitely many constants 
       such that $\scrD := \Pol(\mC)$ 
       satisfies either 1.\ or 2.: 
    \begin{enumerate}
        \item there is a continuous clone homomorphism $\scrD \to\Projs$;
        \item $\scrD$ contains a cyclic 
        (equivalently: a Siggers,
         or a weak near-unanimity) operation $f$
        modulo unary operations of $\scrC$; moreover, $f$ is
        canonical with respect to  {$\Aut(\mC)$}.
     \end{enumerate}
\end{cor}
\begin{proof} 
Let $U_1,\dots,U_n$ be a partition of $A$ such that $\mA$ is a reduct of $(A;U_1,\dots,U_n)$.
If the model-complete core $\mB$ of $\mA$ is finite, then we can expand by a constant for each element of $\mB$,
and the statement follows from Theorem~\ref{thm:finite}.
Otherwise, for some stabilised partition $\{V_1,\dots,V_m\}$ of $B$ the structure $\mB$ is a reduct of $(B;V_1,\dots,V_m)$, by Lemma~\ref{lem:mc-core}.
Then by Proposition~\ref{prop:stabilising-cores}, there are finitely many constants $c_1,\dots,c_m$ such that \purple{$(\mB,c_1,\dots,c_m)$ satisfies the
hypothesis of Theorem~\ref{thm:ua-result-restated}, and the statement follows directly from Theorem~\ref{thm:ua-result-restated}.}
\ignore{the unary part of $\scrD := \Pol(\mB,c_1,\dots,c_m)$
is exactly the set of injections preserving the sets $V_1\setminus\{c_1\},\dots,V_m\setminus\{c_m\},\{c_1\},\dots,\{c_m\}$.
The statement follows directly from Theorem~\ref{thm:ua-result-restated} applied to $\scrD$. }
\end{proof}

We are now ready to give a proof of Theorem~\ref{thm:our-tractability}, which implies
Theorem~\ref{thm:main} from the introduction.

\begin{thm-our-tractability}
Let $\mA$ be a finite-signature reduct of a unary structure.
Then $\Csp(\mA)$ is in P if the model-complete core $\mB$ 
    of $\mA$ has a Siggers polymorphism modulo endomorphisms of $\mB$, 
    and is NP-complete otherwise.
\end{thm-our-tractability}
\begin{proof}
	Let $\mA$ be a finite-signature reduct of $(A;U_1,\dots,U_n)$.
    Let $\mB$ be the model-complete core of $\mA$ and let $\mC$ be the expansion of $\mB$ by finitely many constants given by Corollary~\ref{cor:Siggers-iff-no-homo}.
    Since $\mB$ is a model-complete core, the set of automorphisms of $\mB$ is dense in the set of endomorphisms.
    As in the proof of $\ref{itm:pseudo-cyclic}\Rightarrow\ref{itm:no-h1}$ in Theorem~\ref{thm:ua-result-restated}, we can use this fact to prove that
    $\mC$ has a Siggers polymorphism modulo endomorphisms if, and only if, $\mB$ has such a polymorphism.
    \begin{itemize}
    \item 
    If $\mC$ has such a polymorphism, it has a canonical one, by Corollary~\ref{cor:Siggers-iff-no-homo}.
    Let $m\geq 3$ be greater than the arity of any relation of $\mC$.
    Then $\TB[C]{m}(\mC)$ has a Siggers polymorphism, by Lemma~\ref{lem:polymorphisms-type-algebra}.
    The results from~\cite{BulatovFVConjecture,ZhukFVConjecture} imply that the CSP of $\TB[C]{m}(\mC)$ is in P.
    It follows from Theorem~\ref{thm:types-tractable-implies-structure-tractable} that $\Csp(\mC)$
    is in P, too.
    \item If $\mC$ does not have a Siggers polymorphism modulo endomorphisms, then Corollary~\ref{cor:Siggers-iff-no-homo}
    gives a continuous clone homomorphism from $\Pol(\mC)$ to $\Projs$. 
    By Theorem 1 in~\cite{Topo-Birk}, there exists a polynomial-time reduction
    from, say, 3-SAT to $\Csp(\mC)$. Therefore, $\Csp(\mC)$ is NP-complete. \qedhere
    \end{itemize}
\end{proof}

We mention that using the results from~\cite{BPT-decidability-of-definability}, 
it can be shown that the condition in Theorem~\ref{thm:our-tractability} is decidable:
given subsets $U_1,\dots,U_n$ of $\mA$ (given by the sizes of the sets in the boolean algebra they generate), it is easily seen that one can compute
a finite set of bounds for the age of $(A;U_1,\dots,U_n)$.
Given first-order formulas that define the relations of $\mA$ over $(A;U_1,\dots,U_n)$, it is also possible to compute the model-complete core $\mB$ of $\mA$.
Our results then imply that $\mB$ has a Siggers polymorphism modulo endomorphisms if, and only if, it has a canonical one.
Testing the existence of a canonical function is then decidable, using the results from~\cite{BPT-decidability-of-definability}.

\section{Future Work}
We believe that 
the statement of Corollary~\ref{cor:Siggers-iff-no-homo} also holds for all CSPs expressible in the logic MMSNP introduced by Feder and Vardi~\cite{FederVardi} (MMSNP is a fragment of existential second-order logic). 
Since every MMSNP sentence is equivalent to a finite union of CSPs, this would give another proof that MMSNP has a complexity dichotomy if and only if finite-domain CSPs have a dichotomy, a result due to Feder and Vardi~\cite{FederVardi,Kun}. The previous proof requires intricate constructions of expander structures~\cite{Kun} and it would be interesting to by-pass this. 

Another exciting open problem is to characterise those MMSNP sentences 
that are equivalent to Datalog programs. Such a characterisation is known 
for finite-domain CSPs~\cite{BoundedWidth}. However, this result 
does not immediately yield the answer to the question for MMSNP since
it is not clear that the reduction of Feder and Vardi~\cite{FederVardi,Kun} preserves
Datalog solvability. On the other hand, MMSNP sentences that
are equivalent to first-order sentences have been characterised recently~\cite{Lutz}.
We conjecture that a CSP in MMSNP is equivalent to a Datalog program
if the model-complete core template $\mA$ of the CSP satisfies the tractability 
condition from Theorem~\ref{thm:abstract-canonical-datalog}, in which case the MMSNP problem can be solved by Datalog,
and that otherwise there is an h1 homomorphism from $\Pol(\mA)$ to
the polymorphism clone of a  {module}, in which case CSP$(\mA)$ cannot be solved
by a Datalog program~\cite{AtseriasBulatovDawar}.

\bibliographystyle{plain}
\bibliography{local}

\def\cprime{$'$} \def\cprime{$'$}
\begin{thebibliography}{10}

\bibitem{AtseriasBulatovDawar}
Albert Atserias, Andrei~A. Bulatov, and Anuj Dawar.
\newblock Affine systems of equations and counting infinitary logic.
\newblock {\em Theoretical Computer Science}, 410(18):1666--1683, 2009.

\bibitem{BKOPP}
Libor Barto, Michael Kompatscher, Miroslav Ol\v{s}\'{a}k, Michael Pinsker, and
  Trung~Van Pham.
\newblock {The Two Dichotomy Conjectures for Infinite Domain Constraint
  Satisfaction Problems Are Equivalent}.
\newblock In {\em {Proceedings of the 32nd Annual IEEE Symposium on Logic in
  Computer Science (LICS'17)}}, 2017.

\bibitem{BoundedWidth}
Libor Barto and Marcin Kozik.
\newblock Constraint satisfaction problems of bounded width.
\newblock In {\em Proceedings of the 50th Annual Symposium on Foundations of
  Computer Science (FOCS'09)}, pages 595--603, 2009.

\bibitem{Cyclic}
Libor Barto and Marcin Kozik.
\newblock Absorbing subalgebras, cyclic terms and the constraint satisfaction
  problem.
\newblock {\em Logical Methods in Computer Science}, 8/1(07):1--26, 2012.

\bibitem{wonderland}
Libor Barto, Jakub Opr\v{s}al, and Michael Pinsker.
\newblock The wonderland of reflections.
\newblock Preprint arXiv:1510.04521, 2015.

\bibitem{BartoPinskerDichotomy}
Libor Barto and Michael Pinsker.
\newblock The algebraic dichotomy conjecture for infinite domain constraint
  satisfaction problems.
\newblock In {\em Proceedings of the 31st Annual IEEE Symposium on Logic in
  Computer Science (LICS'16)}, pages 615--622, 2016.
\newblock Preprint arXiv:1602.04353.

\bibitem{Bir-On-the-structure}
Garrett Birkhoff.
\newblock On the structure of abstract algebras.
\newblock {\em Mathematical Proceedings of the Cambridge Philosophical
  Society}, 31(4):433--454, 1935.

\bibitem{Cores-journal}
Manuel Bodirsky.
\newblock Cores of countably categorical structures.
\newblock {\em Logical Methods in Computer Science}, 3(1):1--16, 2007.

\bibitem{Bodirsky-HDR}
Manuel Bodirsky.
\newblock Complexity classification in infinite-domain constraint satisfaction.
\newblock M\'emoire d'habilitation \`a diriger des recherches, Universit\'{e}
  Diderot -- Paris 7. Available at arXiv:1201.0856, 2012.

\bibitem{BodDalJournal}
Manuel Bodirsky and V\'ictor Dalmau.
\newblock Datalog and constraint satisfaction with infinite templates.
\newblock {\em Journal on Computer and System Sciences}, 79:79--100, 2013.
\newblock A preliminary version appeared in the proceedings of the Symposium on
  Theoretical Aspects of Computer Science (STACS'05).

\bibitem{Phylo-Complexity}
Manuel Bodirsky, Peter Jonsson, and Trung~Van Pham.
\newblock The complexity of phylogeny constraint satisfaction.
\newblock In {\em Proceedings of the Symposium on Theoretical Aspects of
  Computer Science (STACS)}, 2016.
\newblock Preprint arXiv:1503.07310.

\bibitem{ecsps}
Manuel Bodirsky and Jan K\'ara.
\newblock The complexity of equality constraint languages.
\newblock {\em Theory of Computing Systems}, 3(2):136--158, 2008.
\newblock A conference version appeared in the proceedings of Computer Science
  Russia {(CSR'06)}.

\bibitem{tcsps-journal}
Manuel Bodirsky and Jan K\'ara.
\newblock The complexity of temporal constraint satisfaction problems.
\newblock {\em Journal of the ACM}, 57(2):1--41, 2009.
\newblock An extended abstract appeared in the Proceedings of the Symposium on
  Theory of Computing (STOC).

\bibitem{BodMarMot}
Manuel Bodirsky, Barnaby Martin, and Antoine Mottet.
\newblock Constraint satisfaction problems over the integers with successor.
\newblock In {\em Proceedings of ICALP}, pages 256--267, 2015.
\newblock ArXiv:1503.08572.

\bibitem{BodMarMot-temporal}
Manuel Bodirsky, Barnaby Martin, and Antoine Mottet.
\newblock {Discrete Temporal Constraint Satisfaction Problems}.
\newblock To appear in the Journal of the ACM, 2017.

\bibitem{BodPin-Schaefer-both}
Manuel Bodirsky and Michael Pinsker.
\newblock Schaefer's theorem for graphs.
\newblock {\em Journal of the ACM}, 62(3):52 pages (article number 19), 2015.
\newblock A conference version appeared in the Proceedings of STOC 2011, pages
  655--664.

\bibitem{Topo-Birk}
Manuel Bodirsky and Michael Pinsker.
\newblock Topological {B}irkhoff.
\newblock {\em Transactions of the American Mathematical Society},
  367:2527--2549, 2015.

\bibitem{canonical}
Manuel Bodirsky and Michael Pinsker.
\newblock Canonical functions: a proof via topological dynamics.
\newblock Preprint arXiv:1610.09660, 2016.

\bibitem{BPP-projective-homomorphisms}
Manuel Bodirsky, Michael Pinsker, and Andr\'{a}s Pongr\'acz.
\newblock Projective clone homomorphisms.
\newblock Preprint arXiv:1409.4601, 2014.

\bibitem{BPT-decidability-of-definability}
Manuel Bodirsky, Michael Pinsker, and Todor Tsankov.
\newblock Decidability of definability.
\newblock {\em Journal of Symbolic Logic}, 78(4):1036--1054, 2013.
\newblock A conference version appeared in the Proceedings of LICS 2011.

\bibitem{equiv-csps}
Manuel Bodirsky and Micha\l Wrona.
\newblock Equivalence constraint satisfaction problems.
\newblock In {\em Proceedings of Computer Science Logic}, volume~16 of {\em
  LIPICS}, pages 122--136. Dagstuhl Publishing, September 2012.

\bibitem{DBLP:journals/corr/BojanczykKL14}
Miko{\l}aj Boja{\'{n}}czyk, Bartek Klin, and Slawomir Lasota.
\newblock Automata theory in nominal sets.
\newblock {\em Logical Methods in Computer Science}, 10(3), 2014.

\bibitem{DBLP:conf/lics/BojanczykKLT13}
Miko{\l}aj Boja{\'{n}}czyk, Bartek Klin, Slawomir Lasota, and Szymon
  Toru{\'{n}}czyk.
\newblock Turing machines with atoms.
\newblock In {\em 28th Annual {ACM/IEEE} Symposium on Logic in Computer
  Science, {LICS} 2013, New Orleans, LA, USA.}, pages 183--192, 2013.

\bibitem{BroxvallJonsson}
Mathias Broxvall and Peter Jonsson.
\newblock Point algebras for temporal reasoning: Algorithms and complexity.
\newblock {\em Artificial Intelligence}, 149(2):179--220, 2003.

\bibitem{BulatovFVConjecture}
Andrei~A. Bulatov.
\newblock {A dichotomy theorem for nonuniform CSPs}.
\newblock In {\em Proceedings of FOCS'17}, 2017.
\newblock {arXiv:1703.03021}.

\bibitem{BulatovJeavons}
Andrei~A. Bulatov and Peter Jeavons.
\newblock Algebraic structures in combinatorial problems.
\newblock Technical report MATH-AL-4-2001, Technische Universit\"at Dresden,
  2001.

\bibitem{JBK}
Andrei~A. Bulatov, Andrei~A. Krokhin, and Peter~G. Jeavons.
\newblock Classifying the complexity of constraints using finite algebras.
\newblock {\em SIAM Journal on Computing}, 34:720--742, 2005.

\bibitem{FederVardi}
Tom\'as Feder and Moshe~Y. Vardi.
\newblock The computational structure of monotone monadic {SNP} and constraint
  satisfaction: {a} study through {D}atalog and group theory.
\newblock {\em {SIAM} Journal on Computing}, 28:57--104, 1999.

\bibitem{Lutz}
Christina Feier, Carsten Lutz, and Antti Kuusisto.
\newblock {Rewritability in Monadic Disjunctive Datalog, MMSNP, and Expressive
  Description Logics.}
\newblock In {\em Proceedings of the 20th International Conference on Database
  Theory (ICDT17)}, 2017.

\bibitem{Hodges}
Wilfrid Hodges.
\newblock {\em A shorter model theory}.
\newblock Cambridge University Press, Cambridge, 1997.

\bibitem{JeavonsClosure}
Peter Jeavons, David Cohen, and Marc Gyssens.
\newblock Closure properties of constraints.
\newblock {\em Journal of the ACM}, 44(4):527--548, 1997.

\bibitem{RCC5JD}
Peter Jonsson and Thomas Drakengren.
\newblock A complete classification of tractability in {RCC}-5.
\newblock {\em Journal of Artificial Intelligence Research}, 6:211--221, 1997.

\bibitem{JonssonThapper15}
Peter Jonsson and Johan Thapper.
\newblock Constraint satisfaction and semilinear expansions of addition over
  the rationals and the reals.
\newblock {\em CoRR}, abs/1506.00479, 2015.

\bibitem{locally-finite}
Bartek Klin, Eryk Kopczynski, Joanna Ochremiak, and Szymon Toru{\'{n}}czyk.
\newblock Locally finite constraint satisfaction problems.
\newblock In {\em 30th Annual {ACM/IEEE} Symposium on Logic in Computer
  Science, {LICS} 2015, Kyoto, Japan.}, pages 475--486, 2015.

\bibitem{KlinLOT14-short}
Bartek Klin, Slawomir Lasota, Joanna Ochremiak, and Szymon Toru{\'{n}}czyk.
\newblock Turing machines with atoms, constraint satisfaction problems, and
  descriptive complexity.
\newblock In {\em Conference on Computer Science Logic {(CSL)} and Symposium on
  Logic in Computer Science (LICS), {CSL-LICS}'14, Vienna, Austria.}, pages
  58:1--58:10, 2014.

\bibitem{KompatscherPham}
Michael Kompatscher and Trung~Van Pham.
\newblock A complexity dichotomy for poset constraint satisfaction.
\newblock {\em Submitted}, 2016.
\newblock Preprint arxiv:1603.00082.

\bibitem{Maltsev-Cond}
Marcin Kozik, Andrei Krokhin, Matt Valeriote, and Ross Willard.
\newblock Characterizations of several maltsev conditions.
\newblock {\em Algebra universalis}, 73(3):205--224, 2015.

\bibitem{Kun}
G\'abor Kun.
\newblock Constraints, {MMSNP}, and expander relational structures.
\newblock {\em Combinatorica}, 33(3):335--347, 2013.

\bibitem{MarotiMcKenzie}
Mikl\'os Mar\'oti and Ralph McKenzie.
\newblock Existence theorems for weakly symmetric operations.
\newblock {\em Algebra Universalis}, 59(3), 2008.

\bibitem{PestovTheorem}
Vladimir Pestov.
\newblock On free actions, minimal flows and a problem by {E}llis.
\newblock {\em Transations of the {A}merican {M}athematical {S}ociety},
  350(10):4149--4165, 1998.

\bibitem{Siggers}
Mark~H. Siggers.
\newblock A strong {M}al'cev condition for varieties omitting the unary type.
\newblock {\em Algebra Universalis}, 64(1):15--20, 2010.

\bibitem{ZhukFVConjecture}
Dmitriy Zhuk.
\newblock {The Proof of CSP Dichotomy Conjecture}.
\newblock In {\em Proceedings of FOCS'17}, 2017.
\newblock {arXiv:1704.01914}.

\end{thebibliography}

\end{document}